\documentclass[review,10pt]{elsarticle}
\usepackage[utf8]{inputenc}
\usepackage{microtype}
\usepackage[margin=2.5cm]{geometry}
\usepackage{setspace}
\usepackage{siunitx}
\usepackage[pagewise]{lineno}
\usepackage[english]{babel}
\usepackage{color}
\usepackage{booktabs}
\usepackage{longtable}
\usepackage[numbers]{natbib}
\usepackage[most]{tcolorbox} 
\usepackage{ragged2e}
\usepackage[overload]{empheq}
\usepackage[justification=justified]{caption}
\usepackage{amssymb}
\newcommand{\Y}{$\checkmark$}
\newcommand{\NN}{$-$}
\usepackage{enumitem}

\usepackage{amsmath,amsfonts,amssymb,mathtools}
\usepackage{cases}
\usepackage{makecell}
\usepackage{multirow}
\usepackage{graphicx,float}
\usepackage{algpseudocode}

\usepackage[]{hyperref}
\hypersetup{colorlinks, linkcolor=blue, citecolor=blue, urlcolor=blue}

\usepackage[ruled,vlined,linesnumbered]{algorithm2e}
\DeclareRobustCommand{\officialeuro}{%
  \ifmmode\expandafter\text\fi
  {\fontencoding{U}\fontfamily{eurosym}\selectfont e}}
\usepackage{framed} % Framing content
\immediate\write18{%
makeindex Main_File.nlo -s nomencl.ist -o Main_File.nls 
}
\usepackage{multicol} % Multiple columns environment
\usepackage{nomencl}[noprefix] % Nomenclature package
\makenomenclature
\renewcommand*\nompreamble{\begin{multicols}{2}}
\renewcommand*\nompostamble{\end{multicols}}

\def\tsc#1{\csdef{#1}{\textsc{\lowercase{#1}}\xspace}}
\tsc{WGM}
\tsc{QE}
\tsc{EP}
\tsc{PMS}
\tsc{BEC}
\tsc{DE}

\begin{document}

    \justifying

    \begin{frontmatter} 
        \title{Forecast-Residual-Based Chance-Constrained Scheduling of Local Energy Communities under PV and Demand Uncertainty}
    
        \author[1]{Franco Cárdenas González}
        \author[1]{Fernando García-Muñoz \corref{mycorrespondingauthor}}
        \cortext[mycorrespondingauthor]{Corresponding author}
        \ead{fernando.garciam@usach.cl}

        \address[1]{University of Santiago of Chile (USACH), Faculty of Engineering, Industrial Engineering Department, Chile}

\begin{abstract}
Day-ahead operation of local energy communities (LECs) is affected by uncertainty in PV generation and electricity demand, which may compromise the feasibility of committed exchanges with the upstream grid. This paper proposes a forecast-residual-based chance-constrained optimization (CCO) framework that characterizes uncertain parameters through their forecasts and associated residual distributions. Under the adopted zero-mean Gaussian assumption, the resulting chance constraints admit an exact deterministic-equivalent reformulation in which uncertainty is embedded through analytical safety margins, preserving the mixed-integer linear structure of the original scheduling problem. The model jointly coordinates PV and BESS operation, low-voltage distribution network constraints, internal energy sharing, and day-ahead grid-exchange commitments. Once the physical schedule is determined, an ex-post allocation stage distributes the available community energy pool among users according to predefined participation coefficients, identifies post-allocation surpluses and deficits, and maximizes their internal matching while preserving the aggregate grid exchanges obtained from the CCO. The framework is evaluated on a reduced 206-node European low-voltage feeder with up to 55 community users and benchmarked against a two-stage stochastic programming formulation. Results show that, for forecast-error levels close to 5\%, the proposed CCO achieves operating outcomes comparable to the stochastic benchmark while retaining computational requirements close to the deterministic formulation. The ex-post allocation results further show that internal matching can reduce gross energy exchanges that require settlement with the upstream grid. %As uncertainty and confidence requirements increase, larger safety margins progressively restrict the feasible operating region.
\end{abstract}
        

\end{frontmatter}

\normalsize

% ============================================================
% SECTION 1 — INTRODUCTION
% ============================================================
\section{Introduction}\label{Intro_Chap}
The sustained decline in the cost of residential-scale distributed energy resources (DERs) has created conditions for widespread prosumer deployment that did not exist a decade ago. The levelized cost of utility-scale solar photovoltaic (PV) generation fell by 90\% between 2010 and 2023 according to IRENA \cite{IRENA2024costs}, while lithium-ion battery energy storage system (BESS) costs declined by 93\% over the same period, reaching approximately USD 197/kWh by 2024 \cite{IRENA2024costs}. These reductions have made joint PV generation and battery storage economically viable for a growing share of residential consumers, who increasingly generate, store, and self-consume electricity on-site. As residential DER deployment scales up, the traditional consumer role is progressively shifting toward active participation in coordinated community operation. Local energy communities (LECs) have emerged as the natural organizational form for this distributed infrastructure, coordinating generation, storage, and consumption at the neighborhood level \cite{zhang2017review}. Within these communities, energy-sharing mechanisms allow locally available energy to be allocated among members, increasing the internal utilization of distributed resources, supporting community-level self-sufficiency, and reducing dependence on exchanges with the upstream grid \cite{tushar2020peer}.

In this context, coordinating LEC operations within broader electricity market frameworks introduces challenges absent in traditional passive-consumer settings. In particular, integrating LECs into upstream electricity markets requires advance commitments of the aggregate energy to be exchanged with the main grid. These commitments provide the information required by electricity suppliers and distribution system operators (DSOs) to carry out their own operational planning and network management. However, because future PV generation and electricity demand are inherently uncertain \cite{soto2021peer}, the quantities committed ahead of operation may deviate from the actual community power balance. This uncertainty may compromise the feasibility and economic efficiency of the committed exchange schedules, potentially leading to network constraint violations or additional balancing market actions.

Addressing this uncertainty has motivated a broad range of optimization methodologies, each relying on different representations of uncertain PV generation and electricity demand to determine robust or risk-aware energy commitments. These formulations typically characterize uncertainty through scenarios \cite{garcia2023benders, nguyen2024advanced}, uncertainty sets \cite{khodoomi2023robust, chang2022vertex}, or probability distributions \cite{aghamohammadi2024decentralized, Li2021ddDRO} describing the possible realizations of the underlying variables. These uncertainty representations are particularly suitable during the early stages of LEC deployment, when limited historical observations are available to statistically characterize the uncertainty affecting community operation. However, as LECs become operational, their routine day-ahead scheduling generates an increasingly rich historical database comprising both the forecasts used for market commitments and their corresponding realizations. This continuously expanding information provides an opportunity to characterize uncertainty directly from the empirical behavior of forecast residuals, rather than relying exclusively on scenario-based representations. In this setting, chance-constrained optimization (CCO) provides an alternative framework for exploiting this information by incorporating the empirical distribution of forecast residuals into the operational scheduling problem. Moreover, when the forecast residuals satisfy suitable distributional assumptions, the chance-constrained formulation admits an exact deterministic-equivalent reformulation that is well suited to large-scale community scheduling.

Nevertheless, despite these advantages, the literature has paid limited attention to applying CCO to operational LEC scheduling while jointly accounting for DER coordination, distribution network (DN) constraints, and community energy sharing. This paper addresses this gap by proposing a chance-constrained day-ahead scheduling framework that models PV generation and electricity demand uncertainty through forecast residuals. By assuming Gaussian forecast residuals, the resulting chance-constrained formulation admits an exact deterministic-equivalent mixed-integer linear reformulation in which forecast uncertainty is analytically embedded through standard-deviation terms that propagate its effects to both local operational variables and network constraints. The resulting physical schedule is subsequently used in an ex-post energy-sharing stage that allocates the available community resources among participants while preserving the aggregate grid exchanges determined by the CCO. Consequently, community operators can exploit historical forecast information through a computationally efficient deterministic formulation without resorting to scenario-based representations of uncertainty.

The remainder of this paper is organized as follows. Section~\ref{sec:literature} reviews the literature on uncertainty-aware operational scheduling of LECs and identifies the corresponding research gaps and contributions. Section~\ref{sec:local_Community_model} presents the deterministic LEC scheduling model and develops the proposed CCO formulation together with its deterministic-equivalent reformulation, uncertainty propagation model, and ex-post energy-sharing allocation mechanism. Section~\ref{sec:case_study} describes the test system, benchmark methodology, and numerical results. Finally, Section~\ref{sec:conclusions} summarizes the main findings, discusses the framework's limitations, and outlines future research directions.

% ============================================================
% SECTION 2 — LITERATURE REVIEW
% ============================================================
\section{Literature Review}\label{sec:literature}
Energy sharing within LECs has been investigated through different organizational and market-based mechanisms aimed at coordinating DERs among participating users. Among these approaches, peer-to-peer (P2P) energy trading has received particular attention as a mechanism in which individual participants exchange energy directly or through a local market. Early research focused primarily on market design, seeking to determine energy prices and trading quantities while ensuring economically efficient interactions among community participants. In this context, researchers proposed different market-clearing mechanisms based on Leader--Follower games \cite{malik2024leader}, Nash bargaining formulations \cite{mishra2024enhancing}, and auction-based schemes such as continuous double auctions \cite{rowe2023optimal}. These approaches mainly aim to maximize either individual profit or community-wide social welfare, often incorporating fairness criteria to ensure an equitable distribution of economic benefits among participants \cite{oh2024peer}. Parallel efforts have investigated different information-exchange architectures, leading to centralized \cite{majidi2019integration}, decentralized \cite{aghamohammadi2024decentralized}, and distributed \cite{wang2023distributed} market implementations, each representing different trade-offs between coordination requirements, computational scalability, and data privacy. Numerous pilot projects worldwide have since demonstrated the practical viability of these market concepts, as documented in recent review papers \cite{zhang2017review}. Nevertheless, because these formulations primarily focused on market transactions, they generally neglected the physical limitations imposed by DNs, leading to economically optimal but, in some cases, operationally infeasible schedules.

Recognizing this limitation, subsequent research extended P2P-based energy-sharing formulations by incorporating explicit DN models into the market-clearing process through optimal power flow (OPF) constraints, thereby jointly optimizing market clearing and network operation \cite{huang2020optimization,zhang2024optimal}. For radial medium- and low-voltage DNs, second-order cone programming (SOCP) relaxations rapidly became the preferred modeling framework because they provide a favorable compromise between physical accuracy and computational tractability. However, the inclusion of network constraints increases the computational complexity of the resulting optimization problem as the number of community participants grows. Consequently, decomposition algorithms have become an essential component of large-scale P2P market formulations, where ADMM is the most widely used method for decomposing the market-clearing problem into agent-level subproblems to preserve data privacy \cite{mohiti2019decentralized,khojasteh2023distributed}. More recently, distributed formulations based on Benders decomposition have been proposed to coordinate both trading decisions and network operation while maintaining a degree of decentralization \cite{garcia2023benders}. Although these decomposition methodologies enable distributed market clearing while preserving the privacy of local information, the results reported in these studies show that the efficient operation of P2P markets relies on sufficient flexibility to compensate for forecast deviations and alleviate network congestion.

In this context, BESS have become a primary source of operational flexibility in LECs, enabling energy arbitrage, reshaping power exchanges with the upstream grid, alleviating network congestion, increasing local renewable energy utilization, and compensating for forecast errors during system operation \cite{KHOJASTEH2026102287}. Consequently, several optimization frameworks explicitly model storage systems as active community resources that can simultaneously schedule charging, discharging, and energy trading decisions \cite{garcia2022multi}. Additional studies have incorporated battery degradation costs into the market-clearing process to ensure economically sustainable long-term operation \cite{tostado2024risk}, while others have compared individual and shared storage configurations, concluding that community-level storage generally provides higher asset utilization and lower investment costs \cite{shadi2024evaluating,deng2024coordinated}. As BESS became increasingly integrated into energy-sharing and P2P-based community operation, their role as a flexibility resource for absorbing operational uncertainty motivated the incorporation of uncertainty-aware optimization techniques.

Accordingly, researchers have adopted different optimization approaches to explicitly address uncertainty in network-constrained energy-sharing with BESS. Stochastic Programming (SP) constitutes the most widely adopted framework, representing uncertainty through discrete scenarios with associated occurrence probabilities while minimizing the expected operational cost \cite{garcia2023benders,nguyen2024advanced,hou2023peer}. Despite its modeling flexibility, SP formulations often become computationally demanding because solution times grow rapidly with the number of scenarios, making scenario generation and reduction a critical modeling step \cite{nazemi2021uncertainty}. Robust Optimization (RO) addresses this limitation by replacing probabilistic scenarios with deterministic uncertainty sets that guarantee feasibility under worst-case realizations \cite{khodoomi2023robust,chang2022vertex,khojasteh2023distributed,wei2024virtual}. Although this approach provides strong operational guarantees, it generally produces conservative schedules that underutilize the available flexibility resources. Distributionally Robust Optimization (DRO) offers an intermediate alternative by optimizing against families of probability distributions consistent with available statistical information rather than against individual scenarios or deterministic uncertainty sets \cite{aghamohammadi2024decentralized}. Existing DRO formulations have relied on moment-based ambiguity sets and Wasserstein-distance metrics to reduce conservatism while preserving robustness against distributional misspecification \cite{nasiri2023moment,Li2021ddDRO}. 

Although these three paradigms represent the current state of the art in uncertainty-aware energy-sharing optimization, they formulate the optimization problem through explicit uncertainty representations based on scenarios, uncertainty sets, or ambiguity sets. RO and DRO are particularly attractive when historical information is limited, or the underlying probability distribution is only partially characterized, whereas SP continues to rely on scenario generation and reduction even as additional historical observations become available. However, as repeated day-ahead community operations accumulate historical forecasts and their corresponding realizations, and as forecasting techniques improve, uncertainty can be characterized with increasing statistical accuracy, opening the possibility of alternative optimization approaches that exploit this information more effectively.

In this context, CCO has emerged as an attractive alternative for uncertainty-aware operation by directly incorporating probabilistic feasibility requirements into the optimization model. Under suitable distributional assumptions on forecast uncertainty, these constraints admit exact deterministic-equivalent reformulations in which the statistical characterization of uncertainty defines analytical safety margins that guarantee the desired confidence level. Consequently, as operational data accumulate and forecasting models improve, this statistical characterization can be translated directly into deterministic safety margins within the optimization model, avoiding explicit scenario generation while preserving computational tractability. Accordingly, CCO has been applied to different power system operation problems, including BESS economic dispatch \cite{cheng2018chance}, storage capacity planning \cite{yahyasoltani2020chance}, joint energy and reserve market clearing \cite{guo2020chance}, and general DN management under uncertainty \cite{wang2019chance}. Early contributions demonstrated that chance-constrained optimal power flow formulations admit exact deterministic equivalents under Gaussian uncertainty assumptions while preserving convexity \cite{bienstock2014chance}. Subsequent works extended these formulations to voltage regulation, thermal congestion management, and multi-objective coordination of network assets and storage systems in active DNs \cite{dall2017chance,zhu2022multi,lee2022chance}. 

However, despite these advantages, applying CCO to energy-sharing problems remains relatively limited. Existing CCO-based contributions have primarily incorporated probabilistic constraints into specific energy-trading or market-clearing formulations, rather than developing comprehensive operational frameworks that jointly capture energy sharing, network operation, DER scheduling, and upstream grid exchanges. Specifically, \cite{guo2020chance} proposed a centralized joint energy-and-reserve market with chance constraints but did not explicitly represent the underlying DN. Similarly, \cite{jia2023peer} formulated probabilistic constraints based on prediction intervals for bilateral multi-microgrid energy trading, whereas \cite{chen2023day} proposed a data-driven chance-constrained framework within a Nash bargaining market architecture. Among the reviewed studies, only \cite{suthar2024chance} combines CCO with network-constrained P2P energy trading, incorporating AC optimal power flow constraints into the market-clearing process. Nevertheless, this formulation neither considers BESS nor addresses the coordinated day-ahead commitment of energy exchanges with the upstream grid, leaving unresolved the joint coordination of storage operation, community energy sharing, DN constraints, and upstream grid commitments that characterizes the day-ahead operation of practical LECs.

Additionally, although most optimization-based energy-sharing formulations reviewed above rely on P2P architectures in which participants retain individual control over their resources and exchange surpluses through bilateral or market-based transactions, current regulatory frameworks for collective self-consumption do not necessarily implement energy sharing through this mechanism. At the European level, the regulatory framework recognizes energy sharing as an activity available to active customers and energy communities \cite{EU2019_944_2024}. In practice, Member States implement this principle through different national allocation and settlement mechanisms. In several established collective self-consumption frameworks, including Spain, Portugal, and France, shared generation is allocated among participating consumers through allocation coefficients rather than through bilateral P2P transactions \cite{MORENO2025e43534}. Under these allocation-based schemes, community generation is assigned to consuming members through predefined or dynamically calculated coefficients, and the allocated energy is reflected in each member's settlement with its retailer. Depending on the national framework, these coefficients may be fixed in advance or calculated after metered consumption and generation become available, as in the more dynamic schemes adopted in Portugal and France. If the allocated energy exceeds a member's consumption, the resulting surplus is handled under applicable national rules, for example by returning it to the originating producer or managing it at the community level. This regulatory structure separates DER and grid operation from subsequent energy accounting, allowing post-allocation surpluses and deficits to be identified and their potential internal matching assessed before settling the residual exchange with the upstream grid.

However, the studies that model energy sharing more closely according to these regulatory allocation mechanisms generally address the allocation problem separately from network-constrained operation under uncertainty. For instance, Mustika et al. \cite{MUSTIKA2022118484} explicitly decouple the physical energy-management stage from the subsequent contractual allocation of community energy, but represent the community as an energy hub without DN constraints and assume deterministic generation and demand profiles. Regulatory-oriented studies have also evaluated alternative allocation coefficients under the Spanish \cite{VILLALONGAPALOU2023101055,GILMENA2023112812} and Portuguese \cite{QUEIROZ2023e20599} frameworks, focusing primarily on the economic performance and efficiency of the resulting sharing schemes. Similarly, Fina \cite{FINA20233822} develops an ex-post allocation algorithm that redistributes information from the DSO-based community allocation according to participants' preferences. These approaches provide practical mechanisms for implementing energy sharing under existing regulatory structures, but they generally treat the allocation stage independently from uncertainty-aware, network-constrained scheduling of the community resources. In particular, none of the reviewed regulation-oriented energy-sharing formulations combines an ex-post allocation mechanism with a CCO-based day-ahead scheduling model that coordinates BESS operation, DN constraints, and upstream grid-exchange commitments.

In summary, the reviewed literature reveals a clear separation between uncertainty-aware energy-sharing formulations, which are predominantly developed around P2P market architectures, and regulation-oriented approaches, which focus on allocation and settlement mechanisms but generally abstract from network-constrained operation under uncertainty. Based on this distinction, four main research gaps emerge. First, existing formulations predominantly rely on SP, RO, or DRO, whereas comparatively little attention has been devoted to CCO despite its ability to explicitly control operational risk while preserving computational tractability. Second, existing CCO-based formulations do not jointly integrate BESS operation, DN constraints, coordinated day-ahead grid-exchange commitments, and an energy-sharing mechanism consistent with allocation-based regulatory frameworks. Third, regulation-oriented energy-sharing studies typically address the allocation and settlement stage separately from network-constrained operation under uncertainty. Finally, computational efficiency remains a key challenge as the number of community participants increases, particularly when uncertainty, DER scheduling, and DN constraints are considered simultaneously. To address these research gaps, this paper makes the following contributions:
 
\begin{enumerate}[label=(\roman*)]
     \item A two-stage LEC operation framework that determines day-ahead upstream-grid exchange commitments while ensuring the probabilistic feasibility of community operation under uncertainty through a CCO formulation based on forecast-residual distributions and a regulation-aligned energy-sharing mechanism. In the first stage, the CCO determines the community’s physical operation under PV-generation and demand uncertainty, while the second stage allocates the resulting community energy among users according to predefined participation coefficients.

     \item A deterministic-equivalent reformulation of the proposed CCO model that analytically characterizes forecast uncertainty, including its cumulative effect on BESS SOC and its spatial propagation through the radial DN, to derive confidence-level safety margins for LEC operation. The resulting physical schedules and available community energy are then used for energy allocation, ensuring that the quantities assigned to users remain consistent with the probabilistically feasible operation determined by the CCO.

     \item A computationally tractable formulation that exploits historical forecast and realization data while preserving the mixed-integer linear structure of the scheduling problem, thereby keeping the computational burden comparatively low relative to scenario-based approaches.
\end{enumerate}

Finally, to facilitate comparison with the existing literature, Table~\ref{tab:litreview} summarizes the main characteristics of all the studies reviewed above. The studies are first grouped according to the energy-sharing mechanism considered. The first group includes works that do not explicitly model community energy sharing; these studies are retained because, given the limited application of CCO to energy-sharing problems, they provide representative examples of CCO formulations developed for related power-system operational problems. The second group comprises P2P-based energy-sharing formulations, whereas the third includes allocation-based approaches using participation coefficients ($\beta$), consistent with the regulatory-oriented mechanisms discussed above. Within each group, the studies are further organized by optimization method, following the sequence deterministic, RO, DRO, SP, TSSP, and CCO. The remaining columns identify the optimization problem type, the operational framework in which it is formulated, the uncertain variables considered and the method used to account for them, whether explicit DN constraints are included, the presence of PV and BESS resources, and the test system used for numerical evaluation.

\begin{table}[h]
\centering
\caption{Comparative summary of the related literature on energy sharing and uncertainty-aware optimization.}
\label{tab:litreview}

\scriptsize
\setlength{\tabcolsep}{1.6pt}
\renewcommand{\arraystretch}{1.10}

\begin{tabular}{c c c c c c c c c c c}
\toprule
\textbf{Ref.} &
\textbf{Year} &
\textbf{E-Sharing} &
\textbf{Problem} &
\textbf{Framework} &
\textbf{Method} &
\textbf{Uncertainty} &
\textbf{N.C.} &
\textbf{PV} &
\textbf{BESS} &
\textbf{Test System}\\
\midrule

\cite{mohiti2019decentralized} & 2019 & -- & CO & DA & RO & Wind, Price & \Y & \NN & \NN & IEEE 33-bus DN \\
\cite{bienstock2014chance} & 2014 & -- & OPF & Intraday & CCO & Wind & \Y & \NN & \NN & 2209--3120-bus TN \\
\cite{dall2017chance} & 2017 & -- & OPF & MPC & CCO & PV, Load & \Y & \Y & \Y & 37-bus DN \\
\cite{cheng2018chance} & 2018 & -- & ED & DA & CCO & Wind & \Y & \NN & \Y & IEEE 24-240-bus TN \\
\cite{wang2019chance} & 2020 & -- & ED & SMPC & CCO & PV, Wind, Load, EV & \Y & \Y & \Y & 12-bus TN \\
\cite{yahyasoltani2020chance} & 2020 & -- & LT-planning & -- & CCO & PV, Wind, Load & \NN & \Y & \Y & MG \\
\cite{nazemi2021uncertainty} & 2021 & -- & Restoration & Post-event & CCO & PV, Wind & \Y & \Y & \Y & IEEE 33/123-bus DN \\
\cite{lee2022chance} & 2023 & -- & OPF & DA & CCO & PV, Load & \Y & \Y & \Y & IEEE 33-bus DN \\
\midrule
\cite{huang2020optimization} & 2020 & P2P & CO & DA & Det. & -- & \NN & \Y & \Y & 30-users MG \\
\cite{garcia2022multi} & 2022 & P2P & CO & DA & Det. & -- & \NN & \Y & \Y & 5-users \\
\cite{goitia2023twostage} & 2023 & P2P & CO & DA/RT & Det. & -- & \NN & \Y & \Y & 15-users \\
\cite{deng2024coordinated} & 2025 & P2P & CO & RT & Det. & -- & \Y & \Y & \Y & IEEE 33-bus DN \\
\cite{wei2024virtual} & 2024 & P2P & CO & DA & Det. & -- & \Y & \Y & \Y & IEEE 13/123-bus DN \\
\cite{tostado2024risk} & 2024 & P2P & CO & DA & RO & PV, Load & \NN & \Y & \Y & 5-users \\
\cite{shadi2024evaluating} & 2025 & P2P & CO & DA & RO & PV, Load & \NN & \Y & \Y & 12-users \\
\cite{khodoomi2023robust} & 2023 & P2P & CO & DA & RO & PV, Wind, Load & \NN & \Y & \Y & 4-users \\
\cite{chang2022vertex} & 2022 & P2P & CO & DA & RO & RES & \Y & \NN & \Y & IEEE 69-bus DN \\
\cite{khojasteh2023distributed} & 2023 & P2P & CO & DA & RO & PV, Price & \NN & \Y & \NN & 12-users \\
\cite{nasiri2023moment} & 2024 & P2P & CO & DA & DRO & PV, Wind & \Y & \Y & \Y & IEEE 33/123-bus DN \\
\cite{aghamohammadi2024decentralized} & 2024 & P2P & CO & DA & DRO & PV & \Y & \Y & \Y & IEEE 37/69-bus DN \\
\cite{Li2021ddDRO} & 2021 & P2P & CO & DA/RT & DRO & PV, Wind, Load & \Y & \Y & \Y & IEEE 123-bus DN \\
\cite{mishra2024enhancing} & 2024 & P2P & CO & DA & SP & PV, Wind & \Y & \Y & \Y & IEEE 33-bus DN \\
\cite{oh2024peer} & 2020 & P2P & CO & DA & SP & PV, Load & \NN & \Y & \NN & 16-users SEC \\
\cite{wang2023distributed} & 2023 & P2P & CO & DA & SP & Wind, Load & \Y & \NN & \Y & IEEE 33-bus DN \\
\cite{hou2023peer} & 2023 & P2P & CO & DA & SP & Wind & \Y & \NN & \Y & IEEE 33-bus DN \\
\cite{garcia2023benders} & 2023 & P2P & CO & DA & TSSP & PV, Load & \Y & \Y & \Y & 69-bus DN \\
\cite{suthar2024chance} & 2024 & P2P & CO & DA & CCO & PV, Load & \Y & \Y & \NN & 15-bus DN \\
\cite{guo2020chance} & 2021 & P2P & RM & DA & CCO & Wind & \Y & \Y & \NN & IEEE 9-bus TN \\
\cite{jia2023peer} & 2023 & P2P & CO & MPC & CCO & PV & \Y & \Y & \Y & IEEE 13/123-bus DN \\
\cite{chen2023day} & 2023 & P2P & CO & DA & CCO & PV, Wind, Load & \Y & \Y & \Y & IEEE 123-bus DN \\
\midrule
\cite{FINA20233822} & 2023 & $\beta$ & EA & Ex-post & Det. & -- & \NN & \Y & \NN & 8-users \\
\cite{GILMENA2023112812} & 2023 & $\beta$ & LT-planning & -- & Det. & -- & \NN & \Y & \Y & 12-users \\
\cite{MUSTIKA2022118484} & 2022 & $\beta$ & EA & Ex-post & Det. & -- & \NN & \Y & \Y & 7-users \\
\cite{VILLALONGAPALOU2023101055} & 2023 & $\beta$ & EA & Ex-post & Det. & -- & \NN & \Y & \NN & 20-users \\
\cite{QUEIROZ2023e20599} & 2023 & $\beta$ & EA & Ex-post & Det. & -- & \NN & \Y & \NN & 4-users \\
\midrule
\textbf{This work} & \textbf{2026} & \textbf{$\beta$} & \textbf{CO + EA} & \textbf{DA} & \textbf{CCO} & \textbf{PV, Load} & \textbf{\Y} & \textbf{\Y} & \textbf{\Y} & \textbf{206-bus LV DN}\\
\bottomrule
\end{tabular}

\vspace{1mm}

\scriptsize
\textbf{Abbreviations.}
$\beta$ = Allocation based on participation coefficients;
N.C. = Network Constraints;
E-Sharing = Energy Sharing;
CO = Community Operation;
EA = Energy Allocation;
OPF = Optimal Power Flow;
ED = Economic Dispatch;
DA = Day-Ahead;
RT = Real-Time;
RM = Reserve market;
MPC = Model Predictive Control;
SMPC = Stochastic Model Predictive Control;
LT-planning = Long-Term Planning;
Det. = Deterministic;
SP = Stochastic Programming;
TSSP = Two-Stage Stochastic Programming;
RO = Robust Optimization;
DRO = Distributionally Robust Optimization;
CCO = Chance-Constrained Optimization;
MG = Microgrid;
DN = Distribution Network;
TN = Transmission Network;
\end{table}

\section{Local energy community model}\label{sec:local_Community_model}
This section presents the proposed two-stage mathematical framework for the day-ahead operation and energy allocation of the LEC. First, we introduce the overall modeling strategy, distinguishing between the community's physical scheduling and the subsequent ex-post allocation of energy among its members. We then present the deterministic formulation as the baseline operational model, followed by its extension to account for PV-generation and demand uncertainty through the proposed CCO approach. Finally, we formulate the ex-post energy-allocation model using the physical schedules and aggregate grid exchanges obtained from the CCO stage. %Then, we discuss its two-stage stochastic programming formulation.

\subsection{Modeling framework}
\begin{figure}[h!]
    \centering
    \includegraphics[width=1.0\linewidth]{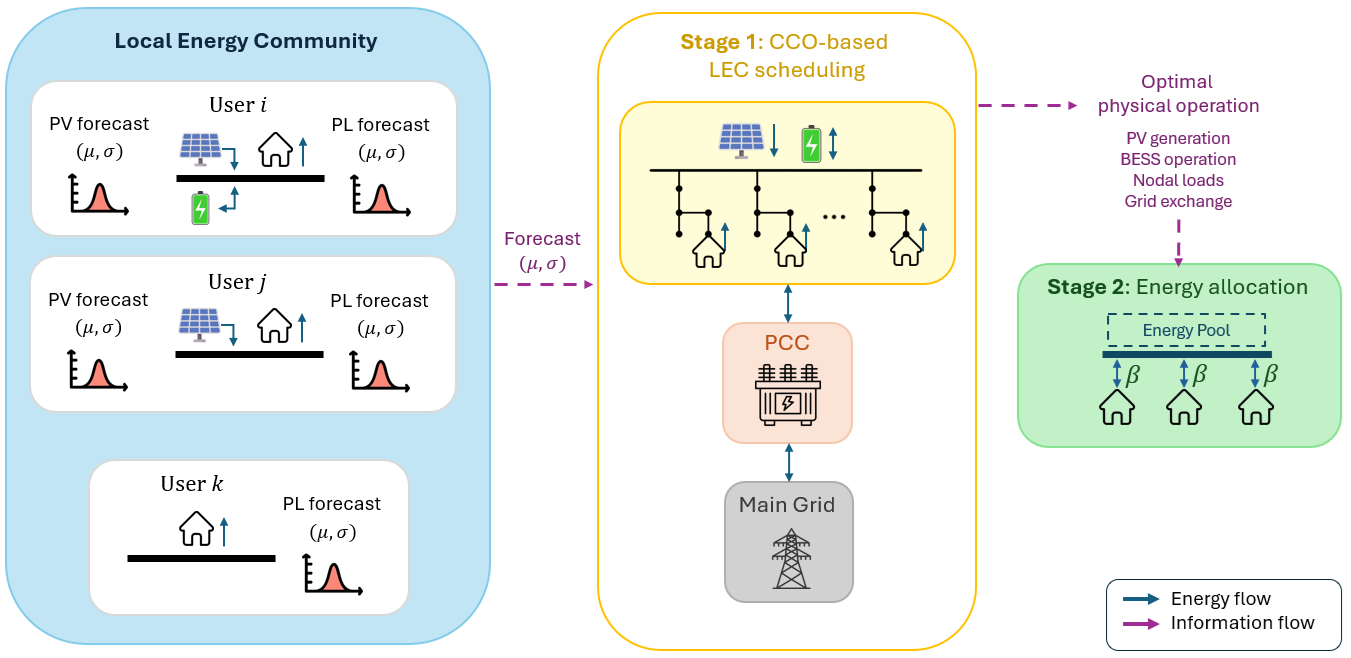}
    \caption{Two-stage framework for CCO-based LEC scheduling and ex-post energy allocation}
    \label{fig:scheme_cco}
\end{figure}

To model the day-ahead operation of the LEC and its energy exchanges with the upstream grid, the proposed framework is organized into two sequential stages, as illustrated in Fig.~\ref{fig:scheme_cco}. This separation reflects the distinction between the community's physical operation and the subsequent accounting-based allocation of energy among its members. Under allocation-based energy-sharing schemes, participation coefficients determine how the community energy is attributed to individual users for settlement and billing purposes, but they do not alter the physical power flows, which are governed by the electrical network and the operation of the available DERs. It is therefore natural to first determine a physically feasible community schedule and the aggregate quantities exchanged with the upstream grid, and subsequently allocate the resulting community energy among users.

Accordingly, the CCO approach adopted in this work exploits historical forecast and realization data to characterize the uncertainty associated with PV generation and demand through their forecast residuals. This information enters the first stage together with the corresponding day-ahead forecasts, allowing the CCO formulation to determine the community's physical operation, including the scheduling of PV and BESS resources, the operation of the radial DN, and the day-ahead grid-exchange commitments while enforcing probabilistic operational feasibility. The resulting physical schedule, including PV generation, BESS operation, nodal demand, and grid exchanges at the points of common coupling (PCC), is then fixed and transferred to the second stage, where the available community energy pool is allocated among users according to predefined participation coefficients and the corresponding post-allocation energy positions are determined.

\subsection{Deterministic model}\label{det_model_LEC}
The LEC is modeled as a low-voltage DN represented by a directed graph $\mathcal{G} = (\mathcal{B}, \mathcal{L})$, where $\mathcal{B}$ denotes the set of buses and $\mathcal{L}$ the set of distribution lines, while a subset $\mathcal{A} \subset \mathcal{B}$ represents the users participating in the energy community. Each agent may operate as a conventional consumer or be equipped with local energy resources, including PV generation systems and combined PV and BESS. The LEC interacts with the upstream grid through one or more PCCs, represented by the subset $\mathcal{E} \subset \mathcal{B}$: each element $i \in \mathcal{E}$ is a network bus at which the community feeder is physically connected to the upstream medium-voltage grid, so the index $i$ in the commitment variables refers to the specific PCC bus at which the exchange takes place. At each of these buses, the community must commit, on a day-ahead basis, the aggregate power to be purchased from and sold to the main grid, denoted by the variables $\kappa_{i,t}^{bg}$ and $\kappa_{i,t}^{sg}$, respectively, where the superscripts $bg$ and $sg$ stand for \emph{buying from} and \emph{selling to} the grid. The model is formulated from the perspective of a local community aggregator that exercises centralized day-ahead coordination on behalf of the LEC. This entity has full information about the technical parameters and forecasts of all agents, jointly determines the dispatch of the community DERs and the corresponding DN operation, and commits the aggregate import/export profile $(\kappa^{bg}_{i,t},\kappa^{sg}_{i,t})$ to the upstream grid.

In this context, the optimal LEC operation is formulated as a Mixed-Integer Linear Programming (MILP) problem that jointly coordinates network operations, DER scheduling, and energy sharing. The scheduling horizon $\mathcal{T}$ comprises the 24 consecutive hourly periods corresponding to the operating day, consistent with the daily clearing cycle of day-ahead electricity markets. The objective is to minimize the total operational cost of the community over the scheduling horizon $\mathcal{T}$, defined as the net cost of energy exchanged with the main grid plus the operational costs associated with local resources, as follows:
\begin{equation} \label{eq:objective}
    \min z = \sum_{t \in \mathcal{T}} \sum_{i \in \mathcal{B}} \left( \lambda_t^{bg} \kappa_{i,t}^{bg} - \lambda_t^{sg} \kappa_{i,t}^{sg} \right) + \sum_{t \in \mathcal{T}} \sum_{i \in \mathcal{A}} C^{op} (pv_{i,t} + ds_{i,t})
\end{equation}

The objective function in \eqref{eq:objective} minimizes the LEC's total operational cost over the scheduling horizon. The first term represents the net cost of energy exchanged with the upstream grid, accounting for the energy purchased from the grid ($\kappa_{i,t}^{bg}$) and the energy sold to the grid ($\kappa_{i,t}^{sg}$), weighted by the corresponding time-dependent buying and selling prices ($\lambda_t^{bg}$ and $\lambda_t^{sg}$). The second term represents the operational cost of local energy resources, with coefficient $C^{op}$ applied to PV generation ($pv_{i,t}$) and battery discharge ($ds_{i,t}$).

We model the network constraints based on the branch-flow formulation in \cite{6507355}, following the radial DN representation introduced in \cite{baran1989network}. This formulation is used to characterize the active- and reactive-power flows, nodal voltage magnitudes, and operating limits of the low-voltage DN as follows:
\begin{subequations} \label{eq:distflow}
\begin{align}
&\sum_{(i,j) \in \mathcal{L}} p_{i,j,t} - \sum_{(j,i) \in \mathcal{L}} p_{j,i,t}
=
\begin{cases}
\Delta p_{i,t}, & \text{if } i \in \mathcal{A}, \\[0.3em]
\kappa_{i,t}^{bg} - \kappa_{i,t}^{sg} & \text{if } i \in \mathcal{E}, \\[0.3em]
0, & \text{Otherwise}.
\end{cases} &&  \forall i \in \mathcal{B}, \forall t \in \mathcal{T} \label{Det_Balance_pg_Eq}\\
&\sum_{(i,j) \in \mathcal{L}} q_{i,j,t} - \sum_{(j,i) \in \mathcal{L}} q_{j,i,t} = qg_{i,t} - QL_{i,t} && \forall i \in \mathcal{B}, \forall t \in \mathcal{T} \label{Det_Balance_qg_Eq}\\
&\Delta p_{i,t} = pv_{i,t} - PL_{i,t} + ds_{i,t} - ch_{i,t}  &&\forall i \in \mathcal{A}, \forall t \in \mathcal{T} \label{Det_DP_Eq}\\
&v_{j,t} = v_{i,t} - 2(R_{i,j} p_{i,j,t} + X_{i,j} q_{i,j,t}) && \forall (i,j) \in \mathcal{L}, \forall t \in \mathcal{T} \label{Det_Voltage_Eq} \\
&A_r p_{i,j,t} + B_r q_{i,j,t} + C_r S_{i,j}^{\max} \le 0  &&\forall (i,j) \in \mathcal{L}, \forall t \in \mathcal{T}, \forall r \in \mathcal{R} \label{Det_AP_Eq}\\
&V_i^{\min} \le v_{i,t} \le V_i^{\max}  && \forall i \in \mathcal{B}, \forall t \in \mathcal{T}\label{Det_LimitV_Eq} \\
&QG_i^{\min} \le qg_{i,t} \le QG_i^{\max} && \forall i \in \mathcal{B}, \forall t \in \mathcal{T} \label{Det_Limit_Q_Eq}\\
&\kappa_{i,t}^{bg} \le PBG^{\max}_i, \quad \kappa_{i,t}^{sg} \le PSG^{\max}_i && \forall i \in \mathcal{E}, \forall t \in \mathcal{T} \label{Det_Limit_Kappas_Eq}
\end{align}
\end{subequations}

Equation (\ref{Det_Balance_pg_Eq}) defines the active power balance at each bus in the network. For user buses, this balance corresponds to the net injected power (see Eq. (\ref{Det_DP_Eq})), which captures the contribution of self-generated PV power, BESS charging and discharging, and user electricity demand. At the PCC with the upstream grid, the balance reflects the net power exchanged with the main grid, i.e., the difference between the aggregated power purchased and the aggregated power sold. For intermediate buses, no net injection is considered. Similarly, Eq. (\ref{Det_Balance_qg_Eq}) establishes the reactive power balance. Following the branch-flow formulation in \cite{6507355} and the radial distribution-network representation introduced in \cite{baran1989network}, we represent the active- and reactive-power balances and the voltage-drop equation under a lossless approximation, neglecting the quadratic loss terms associated with branch currents. Accordingly, Eq. (\ref{Det_Voltage_Eq}) characterizes the voltage drop along each line as a function of the active and reactive power flows and the corresponding line impedance parameters. To ensure secure operation, Eq. (\ref{Det_AP_Eq}) introduces a linearized approximation of the apparent power thermal constraint \cite{yuan2019second}, where $r \in \mathcal{R}$ denotes the set of linear segments used to approximate the circular limit by a polyhedral representation. Furthermore, Eqs. (\ref{Det_LimitV_Eq}--\ref{Det_Limit_Kappas_Eq}) impose the operational limits on voltage magnitude, reactive power capability, and the power exchanged between the energy community and the upstream grid.

Community users may have local DERs, including PV generation and BESS, whose operation determines their net power injection into the network. The following constraints model the technical limits and operational behavior of these prosumer resources over the scheduling horizon.
\begin{subequations}
\begin{align}
    &pv_{i,t} \le PV_t \, \Gamma^{pv}_i &&\forall i \in \mathcal{A}, \forall t \in \mathcal{T} \label{PV_max_Eq} \\
    &soc_{i,t} = soc_{i,t-1} + \phi ch_{i,t} - \frac{1}{\phi} ds_{i,t} &&\forall i \in \mathcal{A}, \forall t \in \mathcal{T} \label{SOC_Eq}\\
    &SOC^{\min} \Gamma^{bt}_i \le soc_{i,t} \le SOC^{\max} \Gamma^{bt}_i &&\forall i \in \mathcal{A}, \forall t \in \mathcal{T} \label{SOC_limit_Eq}\\
    &ch_{i,t} \le PB^{\max} w_{i,t} &&\forall i \in \mathcal{A}, \forall t \in \mathcal{T} \label{BESS_ch_Eq}\\
    &ds_{i,t} \le PB^{\max} (1 - w_{i,t}) - PB^{\max}(1-\nu^{bt}_i) &&\forall i \in \mathcal{A}, \forall t \in \mathcal{T}\label{BESS_ds_Eq}\\    
    &w_{i,t} \le \nu^{bt}_i &&\forall i \in \mathcal{A}, \forall t \in \mathcal{T} \label{BESS_w_Eq}
\end{align}
\end{subequations}
Equation (\ref{PV_max_Eq}) limits the PV power generation ($pv_{i,t}$) at each user to the maximum available capacity ($PV_t$), scaled by the corresponding PV installation parameter ($\Gamma^{pv}_i$). Equation (\ref{SOC_Eq}) describes BESS operation through the state-of-charge balance ($soc_{i,t}$), where $\phi \in (0,1]$ denotes the battery's one-way charge/discharge efficiency. Equation (\ref{SOC_limit_Eq}) enforces the SOC operating range, while Eqs. (\ref{BESS_ch_Eq}) and (\ref{BESS_ds_Eq}) limit the charging ($ch_{i,t}$) and discharging ($ds_{i,t}$) power, respectively. Together with Eq. (\ref{BESS_w_Eq}), these constraints prevent simultaneous charging and discharging through the binary control variable ($w_{i,t}$) and ensure that battery operation is permitted only at nodes equipped with a BESS ($\nu^{bt}_i$).

The net power position of each user is decomposed into surplus and deficit components to characterize how locally available energy and energy requirements are coordinated within the community and with the upstream grid. The following constraints define this decomposition and ensure consistency between user-level exchanges and the aggregate grid commitments.
\begin{subequations}
\begin{align}
    &\Delta p_{i,t} = \Delta p^+_{i,t} - \Delta p^-_{i,t} &&\forall i \in \mathcal{A}, \forall t \in \mathcal{T} \label{Eq_DP}\\
    &\Delta p^+_{i,t} = p^{sg}_{i,t} + p^{sm}_{i,t}&&\forall i \in \mathcal{A}, \forall t \in \mathcal{T} \label{Eq_DP_pos}\\
    &\Delta p^-_{i,t} = p^{bg}_{i,t} + p^{bm}_{i,t}&&\forall i \in \mathcal{A}, \forall t \in \mathcal{T}  \label{Eq_DP_neg}\\
    &\Delta p^+_{i,t} \le \mathcal{M} \, y_{i,t} &&\forall i \in \mathcal{A}, \forall t \in \mathcal{T} \label{Eq_DP_pos_y}\\
    &\Delta p^-_{i,t} \le \mathcal{M} \, (1-y_{i,t})&&\forall i \in \mathcal{A}, \forall t \in \mathcal{T} \label{Eq_DP_neg_y}\\
    &\sum_{i \in \mathcal{A}} p^{sm}_{i,t} = \sum_{i \in \mathcal{A}} p^{bm}_{i,t} \label{Eq_LEM_balance}\\
    &\sum_{i \in \mathcal{A}} p^{bg}_{i,t} = \sum_{i \in \mathcal{E}} \kappa^{bg}_{i,t} \label{Eq_kappa_bg}\\
    &\sum_{i \in \mathcal{A}} p^{sg}_{i,t} = \sum_{i \in \mathcal{E}} \kappa^{sg}_{i,t} \label{Eq_kappa_sg}
\end{align}
\end{subequations}
Equations (\ref{Eq_DP})--(\ref{Eq_DP_neg}) decompose the net power injection ($\Delta p_{i,t}$) of each user into surplus and deficit components. A positive net position ($\Delta p^+_{i,t}$) represents locally available energy that can either be exchanged within the community ($p^{sm}_{i,t}$) or exported to the upstream grid ($p^{sg}_{i,t}$), whereas a negative net position ($\Delta p^-_{i,t}$) represents an energy requirement that can be supplied either through internal community exchanges ($p^{bm}_{i,t}$) or through grid imports ($p^{bg}_{i,t}$). These internal exchange variables maintain consistency with the network power-flow representation and do not constitute the ex-post allocation mechanism described in Section~\ref{sec:energy_sharing}. The binary variable $y_{i,t}$, enforced through Eqs. (\ref{Eq_DP_pos_y}) and (\ref{Eq_DP_neg_y}), ensures that each user has either a surplus or a deficit within a given time period. Equation~(\ref{Eq_LEM_balance}) enforces consistency of the internal energy exchanges by requiring the total energy supplied within the community to equal the total energy received. Finally, Eqs. (\ref{Eq_kappa_bg}) and (\ref{Eq_kappa_sg}) aggregate the user-level grid imports and exports and link them to the corresponding day-ahead power commitments at the PCC.

\subsection{Chance-constrained optimization} \label{subsec:cco}
This section extends the deterministic LEC operational model to account for uncertainty in PV generation and electricity demand through a CCO approach. Under this framework, selected operational limits are enforced probabilistically, requiring them to be satisfied with a prescribed confidence level, while the physical equality constraints remain valid for all uncertainty realizations. The following subsections introduce the general reformulation, derive the corresponding deterministic equivalents, and characterize how uncertainty affects local operational variables and propagates through the DN.

\subsubsection{Problem reformulation}
The deterministic LEC operation model presented in Section~\ref{det_model_LEC} is extended to account for uncertainty in users' PV generation and electricity demand. Let $\tilde{\mathbf{u}}$ denote the vector of uncertain inputs affecting the scheduling problem. Under this representation, the deterministic model can be expressed in the following compact form:

\begin{equation*}
\begin{aligned}
\min_{\mathbf{x}} \quad & \mathbf{c}^{\top}\mathbf{x} \\
\text{s.t.}\quad
& \mathbf{A}\mathbf{x}+\mathbf{E}\tilde{\mathbf{u}}=\mathbf{b},\\
& \mathbf{D}\mathbf{x}+\mathbf{H}\tilde{\mathbf{u}}\leq\mathbf{f},
\end{aligned}
\end{equation*}

where $\mathbf{x}$ collects the decision variables of the deterministic model, while $\tilde{\mathbf{u}}$ collects the realized values of the uncertain inputs, namely users' active power demand ($PL_{i,t}$) and $PV_t$ availability. Matrices $\mathbf{E}$ and $\mathbf{H}$ identify how these uncertain inputs enter the equality and inequality constraints, respectively. The remaining terms retain their usual compact-form interpretation, with $\mathbf{c}$ representing the objective-function coefficients and $\mathbf{A}$ and $\mathbf{D}$ defining the deterministic components of the corresponding constraints.

When historical records of LEC operation are available, including both day-ahead forecasts and the corresponding observed realizations, forecast residuals can be extracted and statistically characterized. In this regard, using these historical forecast and realization data, the uncertainty parameter $\tilde{\mathbf{u}}$ affecting the day-ahead scheduling problem can be modeled as follows:
\begin{equation}\label{uncertainty_eq}
\tilde{\mathbf{u}} = \bar{\mathbf{u}} + \boldsymbol{\xi},
\end{equation}

where $\bar{\mathbf{u}}$ denotes the vector of day-ahead forecasts and $\boldsymbol{\xi}$ is the corresponding forecast-residual vector. The representation in \eqref{uncertainty_eq} follows the standard forecast-error decomposition commonly adopted in time-series forecasting, where the forecast residual is defined as the difference between the realized and forecast values of the uncertain inputs \cite{hyndman2018forecasting}. Under this representation, the uncertainty associated with PV generation and electricity demand is entirely characterized through the statistical properties of the forecast-residual vector, which can be estimated from historical forecast records. These statistical properties are described through the parameters of an underlying probability distribution, such as its moments (e.g., the expected value and variance), which determine the analytical treatment of the associated chance constraints. In particular, some distributions admit exact deterministic-equivalent reformulations, whereas others require approximate reformulations or alternative numerical solution approaches. In this work, the forecast residual associated with each scheduling period $t$ is assumed to follow a Gaussian distribution with zero mean ($\mathbb{E}[\boldsymbol{\xi_t}] = \mathbf{0}$) and period-dependent variance, i.e.,
\begin{equation}
    \xi_t \sim \mathcal{N}(0,\sigma_t^2), \qquad \forall t\in\mathcal{T},
\end{equation}
where the variance $\sigma_t^2$ is estimated from the historical forecast-residual samples corresponding to scheduling period $t$. Note that adopting a Gaussian approximation for the forecast residuals is reasonable in this context, since the residual distribution is estimated separately for each scheduling period rather than over the entire scheduling horizon. Thus, under this assumption, each probabilistic constraint involves a Gaussian random variable whose realization is required to satisfy the operational limits with a prescribed confidence level $(\eta)$, leading to the following chance constraint:
\begin{equation}\label{CCO_Const}
\mathbb{P}\left\{ \mathbf{D}\mathbf{x} + \mathbf{H}\tilde{\mathbf{u}} \le \mathbf{f} \right\}
\ge 1-\eta.
\end{equation}

Since each probabilistic constraint is Gaussian, (\ref{CCO_Const}) admits an exact deterministic-equivalent reformulation \cite{bienstock2014chance} given by
\begin{equation}\label{CCO_Det_Eq_Const}
\mathbf{D}\mathbf{x}
+\mathbf{H}\bar{\mathbf{u}}
+\Phi^{-1}(1-\eta)\boldsymbol{\sigma}
\le
\mathbf{f},
\end{equation}
where $\Phi^{-1}(\cdot)$ denotes the inverse cumulative distribution function of the standard normal distribution. The vector $\boldsymbol{\sigma}$ collects the standard deviations associated with each probabilistic constraint, which are obtained by propagating the forecast-residual uncertainty through the corresponding constraint equations. The analytical expression of each element of $\boldsymbol{\sigma}$ depends on the operational constraint under consideration and is derived in the following subsections. Finally, the resulting deterministic equivalent optimization model is expressed as
\begin{equation}\label{CCO_Det_Eq}
\begin{aligned}
\min_{\mathbf{x}} \quad & \mathbf{c}^\top \mathbf{x} \\
\text{s.t.}\quad
& \mathbf{A}\mathbf{x} + \mathbf{E}\bar{\mathbf{u}} = \mathbf{b},\\
& \mathbf{D}\mathbf{x}
+\mathbf{H}\bar{\mathbf{u}}
+\Phi^{-1}(1-\eta)\boldsymbol{\sigma}
\le
\mathbf{f}.
\end{aligned}
\end{equation}

It is important to note that, in the deterministic equivalent formulation (\ref{CCO_Det_Eq}), the term $\Phi^{-1}(1-\eta)$ only reflects the prescribed risk level and the assumed Gaussian nature of forecast errors, and is therefore independent of the network topology or the location where uncertainty is introduced. The spatial propagation of uncertainty through the DN is captured by the standard deviation vector $\boldsymbol{\sigma}$. While forecast errors originate locally at specific buses due to PV availability and load uncertainty, their impact on network variables such as voltages and line flows accumulates along the network's electrical paths through the power balance and flow equations. Consequently, the value of $\boldsymbol{\sigma}$ differs across constraints, depending on the position of the affected node or line and on the aggregation of downstream uncertainties. Therefore, to extend the deterministic model and explicitly characterize the propagation of uncertainty through the DN, it is necessary to introduce additional structural sets that enable tracking how local forecast errors affect network variables along the radial topology. Thus, let $\mathbf{A}\in\{0,1\}^{|\mathcal{L}|\times|\mathcal{B}|}$ denote a binary line--bus incidence matrix, where $A_{(i,j),k}=1$ if line $(i,j)\in\mathcal{L}$ belongs to the unique electrical path from the substation to bus $k\in\mathcal{B}$, and $A_{(i,j),k}=0$ otherwise. Based on this definition, the set of buses downstream of line $(i,j)$ is defined as
$D(i,j) := \left\{ k \in \mathcal{B} \;\middle|\; A_{(i,j),k} = 1 \right\}$. Similarly, for any bus $i\in\mathcal{B}$, let $\mathcal{P}(i)\subseteq\mathcal{L}$ denote the set of lines belonging to the unique electrical path from the substation to node $i$, defined as $\mathcal{P}(i) := \left\{ (m,n)\in\mathcal{L} \;\middle|\; A_{(m,n),i}=1 \right\}$. These sets allow identifying, respectively, which buses contribute uncertainty to the power flow of a given line and which lines accumulate such uncertainties when computing nodal voltage magnitudes. Through this structure, local forecast errors in active and reactive injections propagate through line flows and subsequently impact network voltage profiles.

\subsubsection{Deterministic equivalent of chance constraints}
Based on the previous discussion, the chance-constrained formulation applies only to operational constraints that may be violated due to forecast uncertainty in load and PV generation. In contrast, equality constraints such as power balance and network flow equations represent physical identities and are therefore enforced deterministically, while their role is to propagate uncertainty across the network. Accordingly, the corresponding probabilistic constraints are expressed as follows:
\begin{subequations}\label{CCO_Prob_set}
\begin{align}
&\mathbb{P}\!\left\{ pv_{i,t} \le PV_t \Gamma^{\mathrm{pv}}_i \right\}
\ge 1-\eta, && \forall i\in\mathcal{A},\ \forall t\in\mathcal{T}, \label{CCO_PV} \\
&\mathbb{P}\!\left\{ pl_{i,t} \ge PL_{i,t} \right\}
\ge 1-\eta, && \forall i\in\mathcal{A},\ \forall t\in\mathcal{T}, \label{CCO_Load} \\
&\mathbb{P}\!\left\{ QG_i^{\min} \le q^{g}_{i,t} \le QG_i^{\max} \right\}
\ge 1-\eta, && \forall i\in\mathcal{B},\ \forall t\in\mathcal{T}, \label{CCO_QG} \\
&\mathbb{P}\!\left\{ v_{i,t} \le V_i^{\max} \right\}
\ge 1-\eta, && \forall i\in\mathcal{B},\ \forall t\in\mathcal{T}, \label{CCO_Volt_max} \\
&\mathbb{P}\!\left\{ v_{i,t} \ge V_i^{\min} \right\}
\ge 1-\eta, && \forall i\in\mathcal{B},\ \forall t\in\mathcal{T}, \label{CCO_Volt_min} \\
&\mathbb{P}\!\left\{ A_r p_{i,j,t} + B_r q_{i,j,t} + C_r S^{\max}_{i,j} \le 0 \right\}
\ge 1-\eta, && \forall (i,j)\in\mathcal{L},\ \forall t\in\mathcal{T},\ \forall r\in\mathcal{R}, \label{CCO_Flow} \\
&\mathbb{P}\!\left\{ SOC^{\min} \Gamma^{bt}_i \le soc_{i,t} \le SOC^{\max} \Gamma^{bt}_i\right\}
\ge 1-\eta, && \forall i\in\mathcal{A},\ \forall t\in\mathcal{T}. \label{CCO_SOC}
\end{align}
\end{subequations}

Note that the auxiliary variable $pl_{i,t}$ is introduced for the purpose of applying the chance-constrained transformation to the load satisfaction constraint. Since the objective function minimizes the operational cost, this variable will always take its lower bound at optimality, i.e., $pl_{i,t}=PL_{i,t}$, and therefore does not alter the solution of the deterministic problem. Thus, applying the deterministic equivalent reformulation in (\ref{CCO_Det_Eq_Const}), the probabilistic constraints in (\ref{CCO_Prob_set}) can be equivalently rewritten as:
\begin{subequations} \label{CCO_Prob_Det_set}
\begin{align}
& pv_{i,t} \;+\; \Phi^{-1}(1-\eta)\;  \Gamma^{pv}_i\sigma^{pv}_{i,t} \;\le\; \bar{PV}_t \,\Gamma^{pv}_i,
&& \forall i\in A,\;\forall t\in T, \\
& pl_{i,t} \;-\; \Phi^{-1}(1-\eta)\;  \sigma^{pl}_{i,t} \;\ge\; \bar{PL}_{i,t} ,
&& \forall i\in A,\;\forall t\in T, \\
& qg_{i,t} \;+\; \Phi^{-1}(1-\eta)\;  \sigma^{qg}_{i,t} \;\le\; QG^{max}_i,
&& \forall i\in B,\;\forall t\in T \\
& qg_{i,t} \;-\; \Phi^{-1}(1-\eta)\;  \sigma^{qg}_{i,t} \;\ge\; QG^{\min}_i,
&& \forall i\in B,\;\forall t\in T \\
& v_{i,t} \;+\; \Phi^{-1}(1-\eta)\;  \sigma^{v}_{i,t} \;\le\; V^{\max}_i,
&& \forall i\in B,\;\forall t\in T, \\
& v_{i,t} \;-\; \Phi^{-1}(1-\eta)\;  \sigma^{v}_{i,t} \;\ge\; V^{\min}_i,
&& \forall i\in B,\;\forall t\in T \\
&A_r\,p_{i,j,t} \;+\;
B_r\,q_{i,j,t} \;+\;
C_r\,S^{\max}_{i,j} +  \Phi^{-1}(1-\eta) \Omega_{i,j,t,r} \;\le\; 0 && \forall (i,j)\in L,\;\forall t \in T,\;\forall r \in R \label{lineal_cco_S}\\
& soc_{i,t} \;+\; \Phi^{-1}(1-\eta)\;  \sigma^{soc}_{i,t} \;\le\;SOC^{\max} \, \Gamma^{bt}_i,
&& \forall i\in A,\;\forall t\in T, \\
& soc_{i,t} \;-\; \Phi^{-1}(1-\eta)\;  \sigma^{soc}_{i,t} \;\ge\;SOC^{\min} \, \Gamma^{bt}_i,
&& \forall i\in A,\;\forall t\in T, \\
% & soc_{i,24} \;-\; \Phi^{-1}(1-\eta)\;  \sigma^{soc}_{i,24} \;\ge\;SOC^{end} \, \Gamma^{bt}_i,
% && \forall i\in A,\;\forall t\in T, \\
&\Omega_{i,j,t,r} = \sqrt{A^2_r \, ({\sigma^{p}}_{i,j,t})^2 + B^2_r \, ({\sigma^{q}}_{i,j,t})^2 +2A_r B_r \mathrm{Cov}(p^{i,j}_t,q^{i,j}_t)} \label{Omega_CCO}
\end{align}
\end{subequations}

Note that each standard deviation term $\sigma$ appearing in \eqref{CCO_Prob_Det_set} is associated with the uncertainty affecting the corresponding constraint and must be computed from the available statistical information. In particular, these terms can be classified into two categories, as previously introduced: i) local uncertainty, which is directly related to
forecast errors at individual buses and is captured by $\sigma^{pv}_{i,t}$, $\sigma^{pl}_{i,t}$, $\sigma^{qg}_{i,t}$, and $\sigma^{soc}_{i,t}$; and ii) propagated uncertainty, which affects network variables such as line power flows and nodal voltages and is represented by $\sigma^{p}_{i,j,t}$ and $\sigma^{v}_{i,t}$, accounting for the accumulation of downstream uncertainties along the radial network topology.

\subsubsection{Local uncertainty}
The $\sigma$ values corresponding to local uncertainty are derived from the forecast errors of PV generation and load demand, denoted by $e^{pv}_{i,t}$ and $e^{pl}_{i,t}$, and are used to compute the corresponding standard deviations as
\begin{equation}\label{sigma_pv_pl}
\sigma^{pv}_{i,t} = \sqrt{\mathrm{Var}\!\left(e^{pv}_{i,t}\right)}, \qquad
\sigma^{pl}_{i,t} = \sqrt{\mathrm{Var}\!\left(e^{pl}_{i,t}\right)},
\end{equation}

At user nodes where two forecast error sources coexist, the overall nodal uncertainty is obtained by combining the previously defined PV and load demand standard deviations \eqref{sigma_delta_p} to capture their combined effect as follows:
\begin{equation}\label{sigma_delta_p}
\sigma^{\Delta p}_{i,t}
=
\sqrt{
\left(\sigma^{pv}_{i,t}\right)^2
+
\left(\sigma^{pl}_{i,t}\right)^2
-
2\,\rho_{pv,pl}\,\sigma^{pv}_{i,t}\sigma^{pl}_{i,t}
},
\end{equation}
Equation~\eqref{sigma_delta_p} shows that the squared terms capture the individual contributions of the PV generation and load demand forecast errors, while the cross term, weighted by the correlation coefficient $\rho_{pv,pl}$, captures their potential statistical dependence. The parameter $\rho_{pv,pl} \in [-1,1]$ represents the linear correlation between PV generation and load demand forecast errors and determines whether their joint effect amplifies or attenuates the resulting nodal uncertainty. In the absence of reliable evidence of a correlation between PV generation and load demand forecast errors, the correlation term can be neglected.

Similarly, the standard deviation associated with the SOC would, in principle, require computing separate standard deviations for the charging power, $\sigma^{ch}_{i,t}$, and the discharging power, $\sigma^{ds}_{i,t}$. However, the binary variable $w_{i,t}$ enforces these variables to operate in mutually exclusive charging and discharging modes. As a result, $ch_{i,t}$ and $ds_{i,t}$ depend on the underlying uncertainty in a non-affine way and are coupled through mutually exclusive operating conditions. This structure leads to chance constraints whose feasible set may be nonconvex and computationally intractable, preventing the derivation of a convex deterministic approximation within the chance-constrained optimization framework \cite{Nemirovski_2007}. 

In this regard, to preserve tractability, we adopt a simplifying approximation based on the dominant operational effect of forecast uncertainty on battery operation. In typical prosumer behavior, positive deviations in PV generation forecasts primarily translate into increased charging, while positive deviations in load demand forecasts primarily translate into increased discharging. Accordingly, PV forecast uncertainty is assumed to mainly affect charging decisions, whereas load demand forecast uncertainty is assumed to mainly affect discharging decisions. Under this approximation, the standard deviation of the charging power is approximated as $\sigma^{ch}_{i,t} \approx \sigma^{pv}_{i,t}$, and the standard deviation of the discharging power is approximated as $\sigma^{ds}_{i,t} \approx \sigma^{pl}_{i,t}$. The resulting uncertainty of the SoC is then obtained by propagating these approximated standard deviations through the SoC dynamic equation, yielding
\begin{equation}\label{sigma_soc}
\sigma^{soc}_{i,t} = \sqrt{ \sum_{\tau=1}^{t} \left[ \phi^2 \left(\alpha^{pv}_{i,\tau}\sigma^{pv}_{i,\tau}\right)^2 + \left(\frac{1}{\phi}\right)^2 \left(\alpha^{pl}_{i,\tau}\sigma^{pl}_{i,\tau}\right)^2 \right] }.
\end{equation}

Note that in Equation~\eqref{sigma_soc}, the index $\tau$ represents the scheduling periods from the beginning of the horizon up to period $t$, thereby accounting for the cumulative effect of forecast uncertainty on the SoC over time. The participation factors $\alpha^{pv}{i,\tau}$ and $\alpha^{pl}{i,\tau}$ represent the fractions of the PV and load forecast errors, respectively, that are absorbed by the BESS and consequently propagated to the SoC. The remaining fractions are assumed to be balanced through other available flexibility mechanisms, such as exchanges with the upstream grid or the local energy exchange. We include no cross-correlation term between PV and load forecast errors because, under the adopted assumption, PV forecast uncertainty affects the SoC through charging decisions, whereas load demand forecast uncertainty affects discharging decisions. Since these effects are modeled as acting through distinct operating modes, including a cross-correlation term would be inconsistent with this approximation. Moreover, the summation assumes temporal independence among forecast residuals across scheduling periods. If either the operating-mode assumption or the temporal-independence assumption were relaxed, the corresponding cross-covariance terms would need to be explicitly incorporated into the expression for $\sigma^{soc}_{i,t}$.

\subsubsection{Propagation of uncertainty} \label{Propagation}
The remaining standard deviation terms correspond to propagated uncertainty resulting from the accumulation of downstream forecast errors along the radial network through the power balance and flow equations. In particular, the deterministic formulation assumes that reactive power required to compensate local reactive demand is generated at the same node through $q^{g}_{i,t}$, as enforced by \eqref{Det_Balance_qg_Eq}. As a result, reactive power compensation carries no uncertainty, and the corresponding standard deviation is $\sigma^{qg}_{i,t}=0$. Consequently, uncertainty propagation through the network is driven exclusively by deviations in net active power injections, which are characterized by the nodal standard deviation $\sigma^{\Delta p}_{i,t}$ defined in \eqref{sigma_delta_p}. Accordingly, $\sigma^{p}_{i,j,t}$ on line $(i,j)$ is obtained by aggregating the net injection uncertainties of all buses electrically downstream of that line. Using the set $\mathcal{D}(i,j)$ previously defined;
\begin{equation}\label{sigma_p}
\sigma^{p}_{m,n,t} = 
\sqrt{
\sum_{k\in\mathcal{D}(m,n)}
\big(\sigma^{\Delta p}_{k,t}\big)^2
+
\sum_{\substack{k\neq l\\k,l\in\mathcal{D}(m,n)}}
\rho_{k,l,t}\,
\sigma^{\Delta p}_{k,t}\,
\sigma^{\Delta p}_{l,t}
}
\qquad \forall (m,n)\in\mathcal{L},\; \forall t
\end{equation}

The second term in \eqref{sigma_p} captures spatial cross-covariances between net injection forecast errors at different downstream buses. However, net injection $\Delta p_{i,t}$ results from combining PV generation and load demand forecast errors and therefore does not directly inherit the spatial correlation structure of PV forecasts alone. In particular, while PV forecast errors may exhibit spatial correlation, load forecast errors are typically node-specific and weakly correlated across prosumers. As a consequence, aggregating PV and load forecast errors into $\Delta p_{i,t}$ reduces the dominance of spatial correlation in the resulting net injection uncertainty. Under this consideration, the cross-covariance terms may be neglected, yielding a diagonal covariance approximation.

The probabilistic line flow constraint in \eqref{lineal_cco_S} enforces the apparent power limit through a linearized combination of active and reactive power flows. Accordingly, the term $\Omega_{i,j,t,r}$ in \eqref{Omega_CCO} represents the standard deviation of the expression $A_r p_{i,j,t} + B_r q_{i,j,t}$. The first two terms in \eqref{Omega_CCO} account for the individual contributions of active and reactive power flow uncertainties, weighted by their respective sensitivity coefficients, while the cross term captures their statistical dependence through the covariance $\mathrm{Cov}(p^{i,j}_t,q^{i,j}_t)$. Since active and reactive power flows on a given line are physically coupled by the power flow equations, this covariance term cannot be neglected a priori and is explicitly retained in the computation of $\Omega_{i,j,t,r}$.

Finally, we derive the standard deviation of the nodal voltage from the deterministic voltage drop equation in \eqref{Det_Voltage_Eq}. In the adopted linearized DN model, voltage deviations are an affine function of active power flows along the feeder, with line resistances as sensitivity coefficients. Consequently, uncertainty in line active power flows propagates linearly to nodal voltages. Using the set $\mathcal{P}(i)$, which denotes the collection of lines along the unique radial path from the slack bus to node $i$, the standard deviation of the voltage magnitude at bus $i$ and time $t$ is computed as
\begin{equation}\label{sigma_v}
\sigma^{v}_{i,t}
=
\sqrt{
\sum_{(m,n)\in \mathcal{P}(i)}
\big(2R_{m,n}\,\sigma^{p}_{m,n,t}\big)^2
},
\qquad \forall i\in \mathcal{B},\; \forall t .
\end{equation}
This expression follows directly from the linear voltage drop relation in \eqref{Det_Voltage_Eq}, where each upstream line contributes proportionally to the nodal voltage deviation. Consequently, voltage uncertainty increases toward downstream buses as additional uncertain line flows accumulate along the radial path.

It is worth noting that the different aggregation paths used to compute $\sigma^{p}_{i,j,t}$ and $\sigma^{v}_{i,t}$ follow directly from the physical structure of power flows in a radial DN. Line power flows are determined by the cumulative net injections of all downstream buses; therefore, uncertainty affecting a given line originates from buses electrically located below that line. In contrast, nodal voltage magnitudes result from the accumulation of voltage drops along the unique path from the slack bus to the node, so voltage uncertainty is driven by uncertainties in all upstream line flows. Accordingly, we use downstream aggregation to compute $\sigma^{p}_{i,j,t}$, whereas upstream aggregation along the feeder path computes $\sigma^{v}_{i,t}$.

With all uncertainty standard deviations thus characterized, the
complete CCO formulation is fully determined. The model minimizes
objective~\eqref{eq:objective} subject to the deterministic equality
constraints~\eqref{eq:distflow}--\eqref{Eq_kappa_sg}, while the
original operational inequalities are replaced by the
deterministic-equivalent bounds~\eqref{CCO_Prob_Det_set}. Each safety margin
$\Phi^{-1}(1{-}\eta)\,\sigma$ is parameterized by the pre-computed
standard deviations $\sigma^{pv}_{i,t}$ and $\sigma^{pl}_{i,t}$ from~\eqref{sigma_pv_pl},
$\sigma^{soc}_{i,t}$ from~\eqref{sigma_soc},
$\sigma^{p}_{i,j,t}$ from~\eqref{sigma_p}, and $\sigma^{v}_{i,t}$
from~\eqref{sigma_v}. The BESS binary logic, grid exchange limits,
and energy sharing arbitrage restrictions retain their original deterministic form.

\subsection{Ex-post energy-sharing allocation} \label{sec:energy_sharing}

Once the CCO scheduling problem is solved, the LEC's physical operation is fixed, and the resulting schedules are transferred to an ex-post energy-sharing stage. Accordingly, the CCO solutions for PV generation $\widehat{pv}_{i,t}$, BESS charging and discharging $\widehat{ch}_{i,t}$ and $\widehat{ds}_{i,t}$, served demand $\widehat{pl}_{i,t}$, and the aggregate grid exchanges $\widehat{\kappa}^{bg}_{i,t}$ and $\widehat{\kappa}^{sg}_{i,t}$ are treated as fixed parameters. This stage therefore does not modify the physical operating point obtained from the CCO, but allocates the available community energy among participating users and determines the corresponding internal and external energy exchanges. The available PV generation and BESS discharge are aggregated into a common energy pool and allocated among community members according to the participation coefficients $\beta_i$. These coefficients account for both the relative contribution of each member to the installed DER capacity and its relative contribution to the community peak demand, as follows:

\begin{equation}
\beta_i=\alpha\beta_i^{\mathrm{DER}}+(1-\alpha)\beta_i^{\mathrm{L}},
\qquad \forall i\in\mathcal{A},
\label{eq:beta_total}
\end{equation}

where

\begin{subequations}
\label{eq:beta_components}
\begin{align}
\beta_i^{\mathrm{DER}} &= \frac{\Gamma_i^{pv}+P^{B\max}\nu_i^{bt}}
{\sum_{j\in\mathcal{A}}\left(\Gamma_j^{pv}+P^{B\max}\nu_j^{bt}\right)},
&& \forall i\in\mathcal{A}, \\
\beta_i^{\mathrm{L}} &= \frac{PL_i^{\max}}
{\sum_{j\in\mathcal{A}}PL_j^{\max}},
&& \forall i\in\mathcal{A}.
\end{align}
\end{subequations}

The parameter $\alpha\in[0,1]$ determines the relative weight assigned to DER ownership and peak demand in the allocation rule. After allocating the common pool according to $\beta_i$, each member obtains a net energy position after accounting for its demand and BESS charging requirements. Members with a positive balance first make their surplus available for sharing within the community, while any remaining surplus is exported to the upstream grid. Conversely, members with a negative balance first cover their deficit through energy shared by other community participants and purchase only the remaining requirement from the grid. Based on this allocation logic, we formulate the ex-post energy-sharing problem as follows.

\begin{equation} \max \; F^{\mathrm{sh}}=\sum_{t\in\mathcal{T}}\sum_{i\in\mathcal{A}} p^{sm}_{i,t} \label{eq:ES_obj} 
\end{equation} 

\begin{subequations} \label{eq:ES_model} 
\begin{align} 
&p^{\mathrm{pool}}_{t} = \sum_{j\in\mathcal{A}}\left(\widehat{pv}_{j,t}+\widehat{ds}_{j,t}\right), && \forall t\in\mathcal{T}, \label{eq:ES_pool} \\ 
&p^{\mathrm{alloc}}_{i,t} = \beta_i p^{\mathrm{pool}}_{t}, && \forall i\in\mathcal{A},\ \forall t\in\mathcal{T}, \label{eq:ES_alloc} \\ 
&\Delta p^{\mathrm{sh}}_{i,t} = p^{\mathrm{alloc}}_{i,t}-\widehat{pl}_{i,t}-\widehat{ch}_{i,t}, && \forall i\in\mathcal{A},\ \forall t\in\mathcal{T}, \label{eq:ES_net_position} \\ 
&\Delta p^{\mathrm{sh}}_{i,t} = \Delta p^{+}_{i,t}-\Delta p^{-}_{i,t}, && \forall i\in\mathcal{A},\ \forall t\in\mathcal{T}, \label{eq:ES_decomposition} \\ 
&\Delta p^{+}_{i,t} = p^{sg}_{i,t}+p^{sm}_{i,t}, && \forall i\in\mathcal{A},\ \forall t\in\mathcal{T}, \label{eq:ES_surplus} \\ 
&\Delta p^{-}_{i,t} = p^{bg}_{i,t}+p^{bm}_{i,t}, && \forall i\in\mathcal{A},\ \forall t\in\mathcal{T}, \label{eq:ES_deficit} \\ 
&\Delta p^{+}_{i,t} \leq M y_{i,t}, && \forall i\in\mathcal{A},\ \forall t\in\mathcal{T}, \label{eq:ES_positive_logic} \\ 
&\Delta p^{-}_{i,t} \leq M(1-y_{i,t}), && \forall i\in\mathcal{A},\ \forall t\in\mathcal{T}, \label{eq:ES_negative_logic} \\ 
&\sum_{i\in\mathcal{A}}p^{sm}_{i,t} = \sum_{i\in\mathcal{A}}p^{bm}_{i,t}, && \forall t\in\mathcal{T}, \label{eq:ES_internal_balance} \\ 
&\sum_{i\in\mathcal{A}}p^{sg}_{i,t} = \sum_{j\in\mathcal{E}}\widehat{\kappa}^{sg}_{j,t}, && \forall t\in\mathcal{T}, \label{eq:ES_grid_export} \\ 
&\sum_{i\in\mathcal{A}}p^{bg}_{i,t} = \sum_{j\in\mathcal{E}}\widehat{\kappa}^{bg}_{j,t}, && \forall t\in\mathcal{T}, \label{eq:ES_grid_import} \\ 
&y_{i,t} \in \{0,1\}, && \forall i\in\mathcal{A},\ \forall t\in\mathcal{T}. \label{eq:ES_binary} 
\end{align} 
\end{subequations}

The objective in (\ref{eq:ES_obj}) maximizes the total energy shared within the community over the scheduling horizon. Since the CCO has already determined the physical dispatch, this stage only reallocates the resulting community resources among the participating members.

Equations~(\ref{eq:ES_pool})--(\ref{eq:ES_net_position}) define the core of the ex-post allocation. First, the fixed PV generation and BESS discharge obtained from the CCO are aggregated into the common energy pool. This pool is then distributed among the community members according to the participation coefficients $\beta_i$. The resulting allocated energy is combined with the fixed served demand and BESS charging schedules to determine the post-allocation net position $\Delta p^{\mathrm{sh}}_{i,t}$ of each member. Equations~(\ref{eq:ES_decomposition})--(\ref{eq:ES_negative_logic}) preserve the same surplus--deficit decomposition and buying/selling logic introduced in the deterministic formulation. Positive post-allocation positions may be assigned either to intra-community sharing or to exports to the upstream grid, while negative positions may be covered through internal sharing or grid purchases. The binary variable $y_{i,t}$ prevents simultaneous surplus and deficit positions. Equation~(\ref{eq:ES_internal_balance}) ensures that total energy sold within the community equals the corresponding internal purchases. Finally, Eqs.~(\ref{eq:ES_grid_export}) and~(\ref{eq:ES_grid_import}) preserve the aggregate grid exchanges obtained from the CCO by enforcing that the sum of member-level exports and imports matches, respectively, the fixed values $\widehat{\kappa}^{sg}_{j,t}$ and $\widehat{\kappa}^{bg}_{j,t}$ at the PCC.

\section{Case Study and Computational Results}
\label{sec:case_study}
This section presents the case study used to evaluate the proposed CCO framework and reports the corresponding computational results. First, we describe the test system, input data, uncertainty calibration procedure, and benchmark formulations. The numerical analysis then examines the proposed model's economic and operational performance under different seasonal conditions, along with the impact of uncertainty propagation on community and network operation. Finally, the sensitivity of the CCO formulation to the forecast-error level and confidence requirement is assessed, with particular attention to its economic performance, feasibility, and computational behavior relative to the reference formulations.

\subsection{Test System and Parameters}
\label{sec:test_system}

\begin{figure}[htbp]
    \centering
    \includegraphics[width=0.80\linewidth]{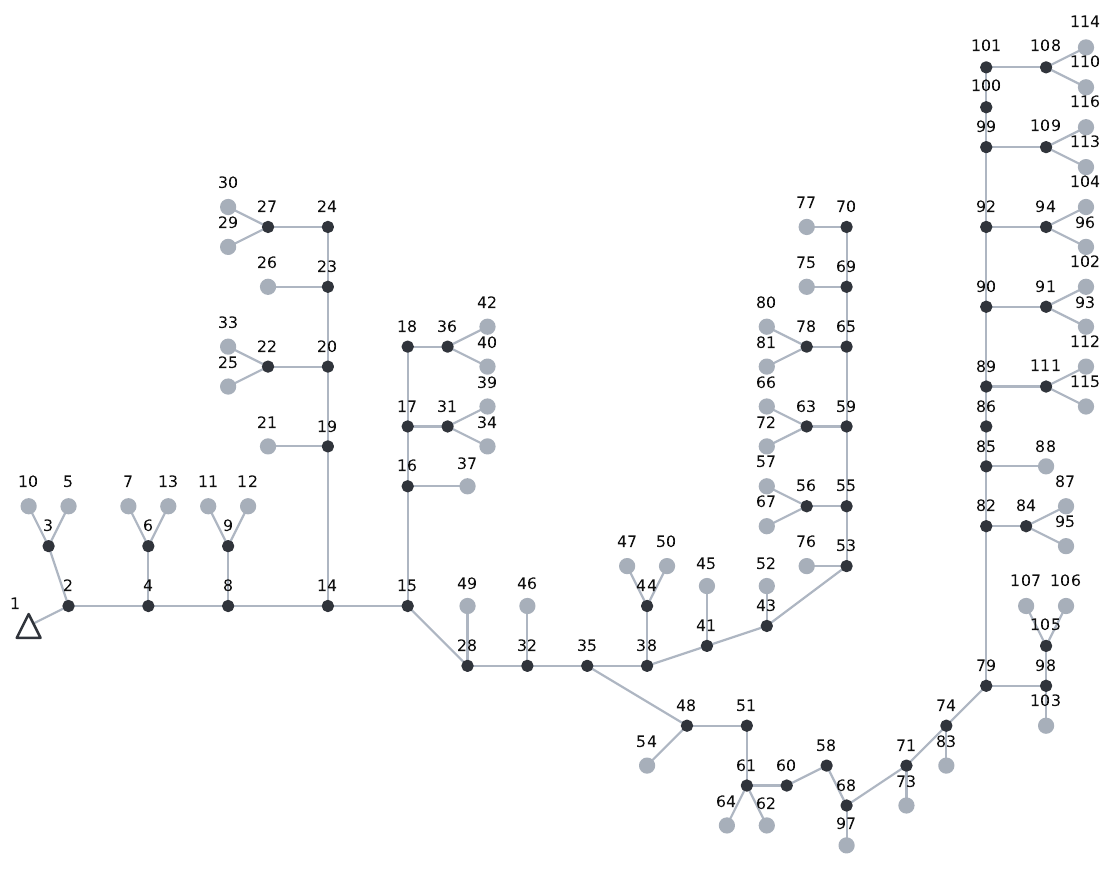}
    \caption{Reduced 206-node system of the IEEE European Low-Voltage Test Feeder.}
    \label{fig:network}
\end{figure}

The network used to evaluate the CCO model presented in the previous section is based on the reduced electrically equivalent model of the IEEE European Low-Voltage Test Feeder proposed in~\cite{Khan2022}, as illustrated in Fig.~\ref{fig:network}.

The network comprises 206 buses connected through 205 distribution branches, with bus 1 representing the point of common coupling with the upstream medium-voltage grid, while the remaining 205 buses host loads and/or DERs. Of these buses, 55 are designated as community agents: 31 are prosumers equipped with rooftop PV systems and BESS, while the remaining 24 are conventional consumers without distributed generation or storage. Table~\ref{tab:parameters} summarizes the network, DER, and operational parameters adopted for the case study.

\begin{table}[htbp]
\centering
\caption{Key parameters of the 206-node test system and distributed energy resources.}
\label{tab:parameters}
\begin{tabular}{@{}l r @{\hspace{1.2cm}} l r@{}}
\toprule
\multicolumn{2}{c}{\textbf{Network Topology}} &
\multicolumn{2}{c}{\textbf{Community Agents}} \\
\cmidrule(r){1-2}\cmidrule(l){3-4}

Total network nodes     & 206 &
Total agents $|\mathcal{A}|$ & 55 \\

Distribution branches   & 205 &
\quad Prosumers (PV+BESS) & 31 \\

PCC bus                  & 1 &
\quad Conventional consumers & 24 \\

Nominal voltage          & \SI{0.4}{kV} &
Peak aggregated demand & $187$~kW \\

Nodal voltage limits     & $[0.90,\;1.1]$~p.u. &
Peak demand per agent & $3.4$~kW \\

                        & &
Installed PV capacity range & $5$--$8$~kW \\

                        & &
Total installed PV capacity & $195.0$~kW \\

\midrule

\multicolumn{2}{c}{\textbf{Battery Energy Storage}} &
\multicolumn{2}{c}{\textbf{Operational and Market Parameters}} \\
\cmidrule(r){1-2}\cmidrule(l){3-4}

BESS capacity range  & $5.0$--$8.8$~kWh &
Scheduling horizon $|\mathcal{T}|$ & 24~h \\

Total installed BESS capacity & $205.5$~kWh &
Time resolution & Hourly \\

BESS efficiency $\phi$ & $0.9$ &
Grid buy price $\lambda^{bg}$ & $11.9$--$24.9$~\textcent/kWh \\

SoC limits & $[10\%,\;90\%]$ &
Feed-in tariff $\lambda^{sg}$ & $6.9$--$18.8$~\textcent/kWh \\

BESS power limit $PB^{\max}$ & $3.4$~kW &
TSSP scenarios $|\mathcal{S}|$ & 9 \\

\bottomrule
\end{tabular}
\end{table}

%CCO confidence levels $(1-\eta)$   & $50,75,90,95,99\%$  \\
%Forecast-error levels           & $5,10,15\%$         \\

To represent the uncertainty considered in the analysis, we generated data randomly. In particular, we generate demand profiles for four seasonal conditions (winter, spring, summer, and autumn) and independently for each user of the energy community using the Python tool in \cite{cadema_household_load_profile}. The tool scales the normalized profiles according to the peak demand assigned to each user in the previously defined network configuration. Because the original residential profiles are generated at a finer temporal resolution, we aggregate them to the hourly resolution used in the optimization model. Similarly, we generate a common daily PV availability profile for the community using \texttt{pvlib} \cite{Anderson2023pvlib}. Thus, for each season, we randomly select a date within the corresponding seasonal interval and classify it as clear, partly cloudy, or cloudy according to season-dependent probabilities. We obtain clear-sky irradiance using the Ineichen model for Madrid and convert it to plane-of-array irradiance for a south-oriented PV installation with a ($30^\circ$) tilt. A stochastic cloud factor is then applied to represent the selected weather condition. The PV profile is initially generated at 15-minute resolution and subsequently averaged to hourly values, resulting in a complete 24-hour seasonal trajectory consistent with the selected day type. Electricity prices for buying and selling are based on historical prices reported in \cite{REE_ESIOS_Prices}. For each simulation and season, we randomly select one complete weekday from the corresponding seasonal period and use its 24-hour buying and selling price profiles in the optimization.

To estimate the PV and user-demand standard deviations required by the CCO model, we used a synthetic calibration procedure. First, we generated 1000 simulations, each containing one 24-hour demand and PV realization for each season, using these profiles as the reference realizations. Second, we constructed synthetic forecasts by introducing zero-mean Gaussian relative perturbations with a prescribed standard deviation (e.g., 5\%, 10\%, or 15\%). After perturbation, we constrained the resulting forecasts to their physically admissible normalized range, particularly to preserve zero PV production during nighttime hours. We then recalculated effective forecast residuals from the final forecast and reference realization values and obtained their sample standard deviations for each season and hour, and, for demand, separately for each community user.

We use these calibrated standard deviations to parameterize the proposed CCO formulation. To evaluate its performance, we compare it against three reference cases: an empirical-($\sigma$) CCO, a two-stage stochastic programming (TSSP) formulation, and a perfect-information (PI) case. Each benchmark serves a different purpose: the PI case provides an ex-post reference under fully revealed uncertainty, the TSSP represents a conventional scenario-based treatment of uncertainty, and the empirical-($\sigma$) CCO assesses the effect of replacing the analytical uncertainty-propagation expressions with empirically estimated standard deviations. It should be noted that the proposed CCO and the TSSP are not equivalent optimization models, since they represent and propagate uncertainty through fundamentally different formulations. The TSSP is therefore included as a benchmark because scenario-based stochastic programming is widely adopted in the related literature, and we use its results only as a reference for comparing objective-function values and computational performance.

For the PI case, we compare the CCO solution with the operation that would have been obtained if actual demand and PV realizations had been known in advance. For each simulation and season, we generate one synthetic realization from the same demand and PV forecasts used by the CCO, using the previously calibrated forecast-error standard deviations. We generate demand realizations independently for each user and hour, while PV realizations follow the corresponding seasonal and hourly uncertainty. We use the resulting profiles directly as deterministic inputs to the optimization model, while the electricity-price profile remains identical to that used by the CCO.

For the TSSP benchmark, we construct the scenario set from the complete demand and PV realizations generated during the 1000-simulation calibration stage. For each simulation and season, we consider the available realizations corresponding to the same seasonal condition and PV day type as the empirical scenario pool. From this pool, we randomly select 200 complete daily realizations with replacement using an empirical Monte Carlo bootstrap procedure. Each selected realization contains the 24-hour demand profiles of all community users together with the corresponding 24-hour PV profile. We then reduce the 200 candidate scenarios to nine representative scenarios using K-means clustering. The centroid of each cluster defines the corresponding representative demand and PV trajectories, and the scenario probability comes from the relative number of candidate realizations assigned to that cluster.

Finally, for the empirical-($\sigma$) CCO, we assess how accurately the analytical expressions used in the proposed CCO reproduce the variability observed in the optimization variables and how any discrepancies affect the objective value and operating decisions. Using the same forecast and realization generation procedure described above, we solve two deterministic problems for each simulation and season: one with the forecast profiles and one with the corresponding realized demand and PV profiles. We collect the resulting differences in net injections, active-power flows, nodal voltages, and battery state of charge across repeated simulations, and use their standard deviations as empirical counterparts to the analytically propagated uncertainty terms. The original demand and PV forecast-error standard deviations remain unchanged, so both CCO variants use the same forecasts, input uncertainty calibration, and confidence level, differing only in the standard deviations used to represent uncertainty propagation through the model.

The complete implementation is publicly available~\cite{CardenasGonzalez_GarciaMunoz_2026}, including the case-study data, the calibrated standard deviations used to parameterize the CCO formulation, and the source code required to reproduce the computational simulations reported in this study. We implemented all models in Python~3.10 using the Pyomo~6.7 framework and solved them with Gurobi Optimizer version~11.0.1. We ran computational experiments on a workstation running Ubuntu~22.04.5~LTS with an AMD Ryzen~7 3700U processor (4 physical cores, 2.2--4.0~GHz, 8 logical threads) and 16~GB DDR4 RAM. 

\subsection{CCO computational performance}\label{Comp_performance}
We run 25 simulations for each season, resulting in 100 simulations across the four seasonal conditions. Table~\ref{tab:benchmark_comparison} reports the average economic and computational performance over these simulations for the proposed CCO formulation, evaluated using a 5\% forecast-error level and a 95\% confidence level, relative to the PI, TSSP, and the alternative CCO formulation in which $\sigma$ is estimated empirically, denoted as CCO*. For this numerical setting, we approximate the spatial cross-covariance terms in Eq.~\eqref{sigma_p} using a diagonal covariance structure, since net-injection uncertainty combines PV and demand forecast errors, with the latter reducing the relative contribution of spatial correlation. 

\begin{table}[hbt]
\centering
\caption{Economic and computational performance relative to the perfect-information benchmark.}
\label{tab:benchmark_comparison}
\begin{tabular}{llcccc}
\toprule
Scope / Season & Model & Mean Obj. & MAE vs. PI & 95\% CI of MAE & Mean time [s] \\
\midrule
\multirow{4}{*}{Overall}
& PI   & -12.4 & --   & --           & 19.4 \\
& CCO  &  -4.0 &  8.3 & [7.3, 9.4]   & 19.7 \\
& TSSP &  -3.7 & 10.0 & [8.1, 12.1]  & 40.1 \\
& CCO* &   2.1 & 14.5 & [12.4, 16.9] & 19.1 \\
\midrule
\multirow{4}{*}{Winter}
& PI   &  2.3 & --  & --          & 20.0 \\
& CCO  &  9.4 & 7.0 & [5.3, 8.9]  & 20.5 \\
& TSSP &  9.3 & 7.9 & [5.1, 11.0] & 34.0 \\
& CCO* & 10.1 & 7.8 & [6.0, 9.5]  & 19.5 \\
\midrule
\multirow{4}{*}{Spring}
& PI   & -6.3 & --   & --           & 19.4 \\
& CCO  & -0.7 &  5.6 & [4.0, 7.3]   & 19.4 \\
& TSSP &  1.8 &  8.7 & [5.6, 12.5]  & 43.3 \\
& CCO* & 21.5 & 27.8 & [22.5, 33.9] & 19.4 \\
\midrule
\multirow{4}{*}{Summer}
& PI   & -37.2 & --   & --          & 19.1 \\
& CCO  & -26.1 & 11.1 & [9.0, 13.2] & 19.4 \\
& TSSP & -34.4 &  6.2 & [4.8, 7.6]  & 41.9 \\
& CCO* & -25.8 & 11.3 & [9.3, 13.4] & 18.7 \\
\midrule
\multirow{4}{*}{Autumn}
& PI   & -8.4 & --   & --           & 19.0 \\
& CCO  &  1.3 &  9.6 & [7.8, 11.6]  & 19.5 \\
& TSSP &  8.6 & 17.2 & [12.0, 23.0] & 41.3 \\
& CCO* &  2.7 & 11.0 & [9.4, 12.9]  & 18.7 \\
\bottomrule
\end{tabular}
\end{table}

In terms of computational burden, both CCO formulations have solution times comparable to the deterministic PI case, requiring an average of 19.7~s for the proposed CCO and 19.1~s for CCO*, compared with 19.4~s for PI. By contrast, the TSSP requires 40.1~s on average, approximately twice the computational time of the proposed CCO. This higher computational burden stems from the TSSP's scenario-based stochastic formulation, which accounts for multiple realizations, whereas the CCO formulations retain a deterministic optimization structure once they incorporate uncertainty margins.

From an economic perspective, the proposed CCO provides the closest overall approximation to the PI benchmark, with an MAE of 8.3, compared with 10.0 for the TSSP and 14.5 for CCO*. However, the seasonal results show that this relative performance varies across operating conditions. The proposed CCO achieves the lowest MAE in winter, spring, and autumn, whereas the TSSP performs better in summer, reducing the MAE to 6.2 compared with 11.1 for the CCO. Conversely, the largest difference occurs in autumn, when the TSSP MAE increases to 17.2 while the CCO remains at 9.6. These seasonal differences suggest that the relative performance of the two uncertainty representations is closely tied to the variability each formulation captures: the dispersion of the generated scenarios in the TSSP and the forecast-error variability represented by $\sigma$ in the CCO, relative to the corresponding PI realization.

\begin{figure}[htbp]
    \centering
    \includegraphics[width=1.0\linewidth]{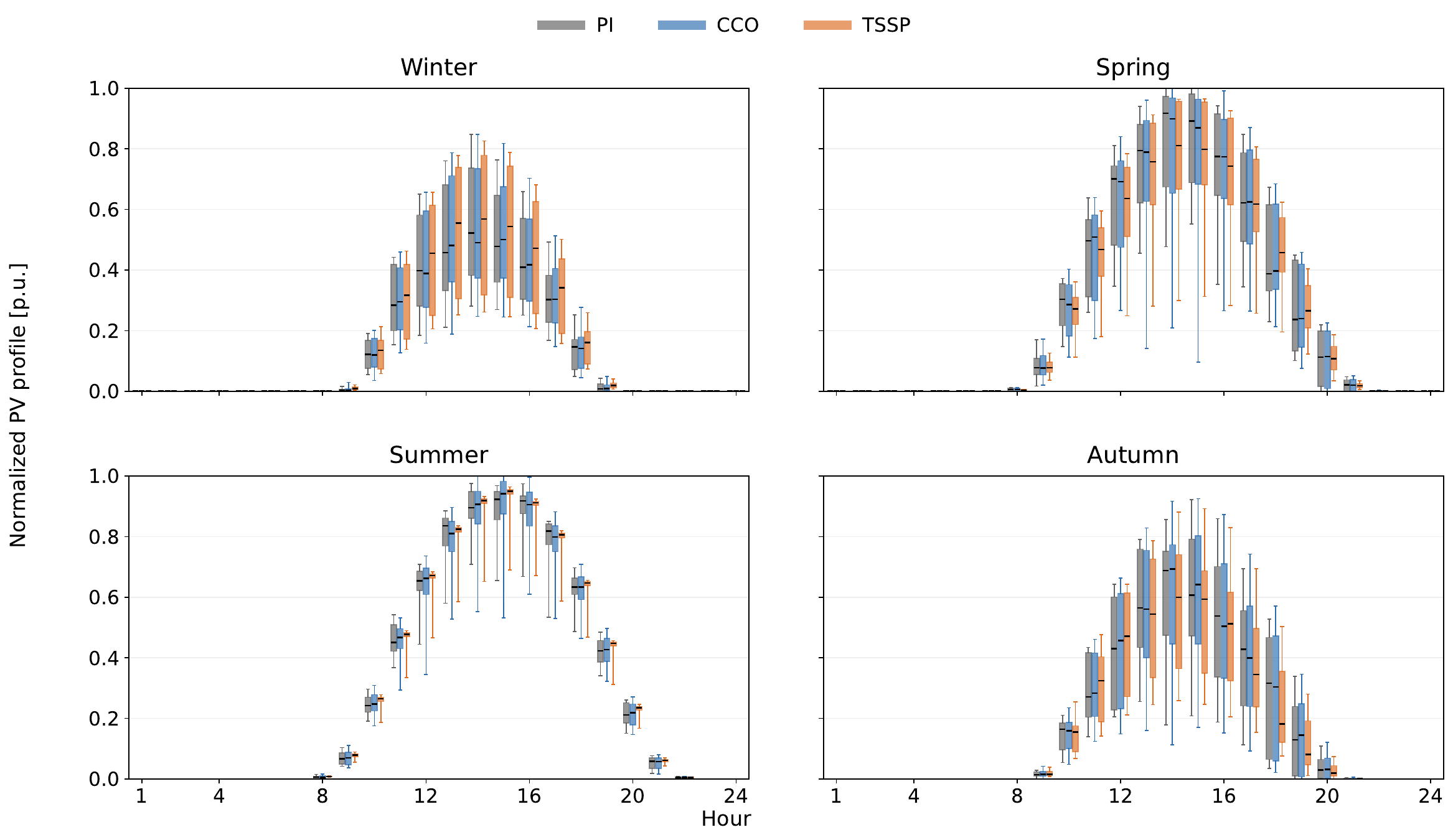}
    \caption{PV uncertainty representation: PI vs. CCO vs. TSSP}
    \label{fig:pv_uncertainty}
\end{figure}

In particular, the TSSP formulation's superior summer performance stems primarily from its representation of PV uncertainty. As shown in Fig.~\ref{fig:pv_uncertainty}, the summer PV scenarios generated for the TSSP exhibit a markedly narrower dispersion than in the other seasons and remain closely aligned with the central tendency of the realized PV profiles throughout the main solar-generation hours. This indicates that, under summer conditions, the scenario set accurately represents realized uncertainty. Consequently, the TSSP provides a more precise description of ex-post PV conditions, which helps explain why its objective value remains closer to the PI benchmark during this season. In contrast, during winter, spring, and autumn, the realized PV distributions exhibit a broader spread, while the CCO uncertainty representation more consistently covers this variability around the forecast. This broader coverage better aligns with the observed dispersion of realized PV profiles in these seasons, consistent with the lower MAE the CCO achieves relative to the TSSP in the corresponding cases.

\begin{figure}[h!]
    \centering
    \includegraphics[width=0.9\linewidth]{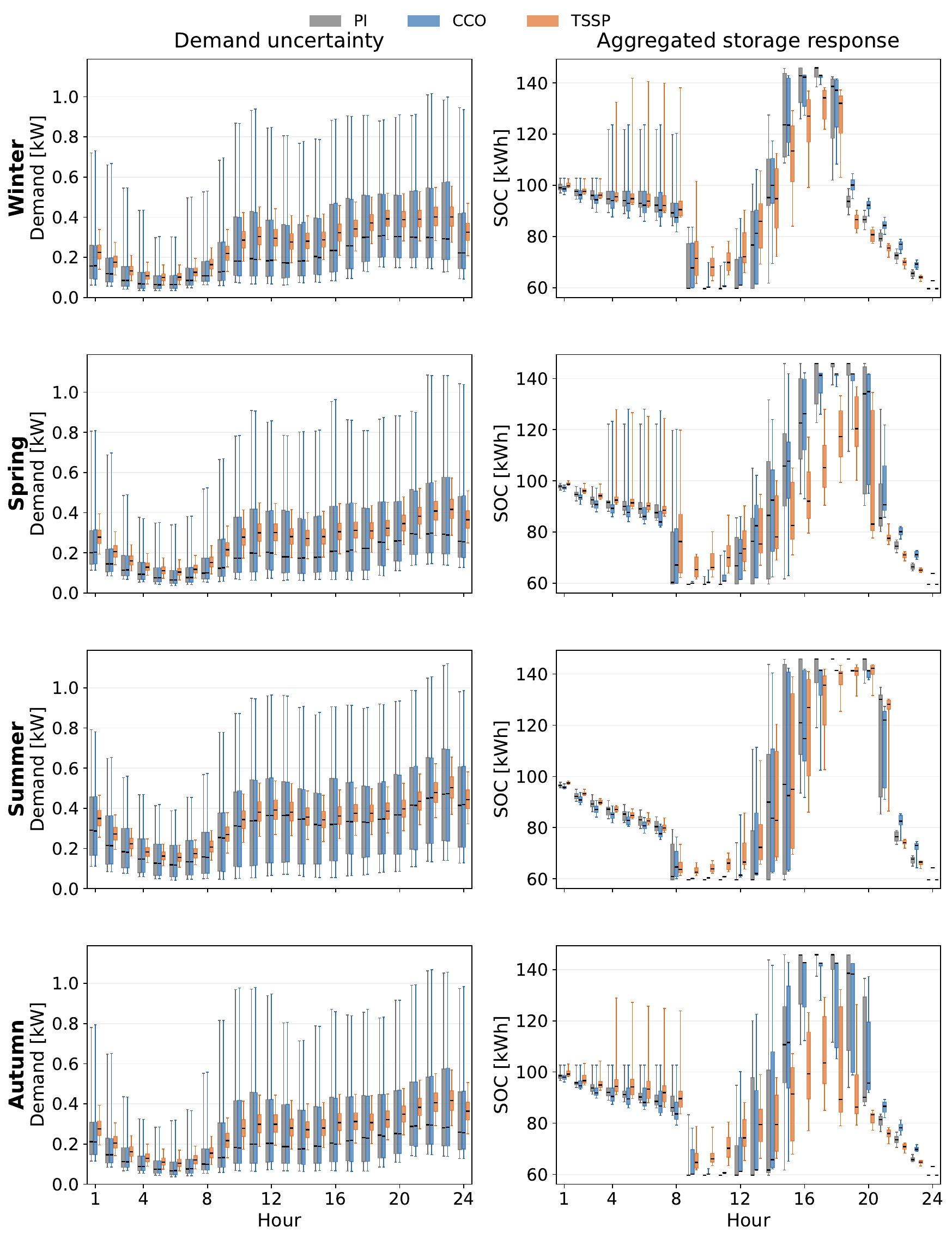}
    \caption{Seasonal demand uncertainty and aggregated storage response}
    \label{fig:demand_soc_uncertainty}
\end{figure}

In contrast, to analyze differences in the autumn season, we need to examine the relationship between the demand representation and the resulting SOC trajectory. Thus, as shown in Fig.~\ref{fig:demand_soc_uncertainty}, the TSSP demand scenarios generally show narrower dispersion than the realized demand distribution and tend to be centered at slightly higher demand levels, particularly during daytime and evening periods. This representation leads the TSSP to rely more heavily on available stored energy, as reflected in the aggregated SOC trajectories shown in Fig.~\ref{fig:demand_soc_uncertainty}. Although the magnitude of this effect varies by season, the TSSP generally maintains a lower SOC than the PI solution, whereas the proposed CCO follows the PI trajectory more closely. Table~\ref{tab:soc_deviations} confirms this behavior quantitatively. The CCO presents SOC MAEs between 3.08 and 4.22~kWh across all seasons, with biases close to zero, whereas the TSSP reaches substantially larger deviations, particularly in spring and autumn, with MAEs of 13.17 and 12.77~kWh and negative biases of $-4.90$ and $-5.78$~kWh, respectively. The economic impact of these deviations is nevertheless season-dependent. In particular, the mean grid buying price is considerably higher in autumn (154.76~EUR/MWh) than in spring (123.71~EUR/MWh). Consequently, the TSSP is penalized in autumn for maintaining a lower SOC, as the associated reduction in community energy autonomy increases exposure to comparatively more expensive grid purchases, which explains the larger deviation of the TSSP objective value from the PI benchmark observed in autumn.

\begin{table}[h]
\centering
\caption{Seasonal deviations of the aggregated community SOC with respect to the perfect-information benchmark [kWh], together with the mean grid buying price [EUR/MWh].}
\label{tab:soc_deviations}
\begin{tabular}{lrrrrr}
\toprule
Season & CCO MAE & TSSP MAE & CCO Bias & TSSP Bias & Mean buying price \\
\midrule
Winter & 3.32 & 8.15  &  1.21 & -0.26 & 149.82 \\
Spring & 4.22 & 13.17 &  1.26 & -4.90 & 123.71 \\
Summer & 3.08 & 5.91  & -0.48 &  0.25 & 151.97 \\
Autumn & 4.04 & 12.77 &  0.23 & -5.78 & 154.76 \\
\bottomrule
\end{tabular}
\end{table}

Additionally, CCO* performs markedly worse in spring mainly because $\sigma_{soc}$ is estimated empirically. In this case, forecast and realized operating conditions change the BESS charging strategy, so differences between the corresponding deterministic SOC trajectories reflect forecast uncertainty and changes in the optimal storage response. Consequently, using these differences to estimate $\sigma_{soc}$ yields larger values than analytical propagation. The resulting SOC safety margins become more restrictive and, in some realizations, approach the BESS energy-capacity limits, forcing the model to adjust its operating decisions to remain feasible. This additional restriction is ultimately reflected in the objective value and explains the pronounced deviation of CCO* observed in spring. This result also indicates that estimating SOC uncertainty directly from differences between optimized trajectories may be unsuitable when forecast deviations trigger changes in the BESS operating strategy.

On the other hand, the DN model adopts a lossless branch-flow approximation and a polyhedral representation of the apparent-power limit to preserve the MILP structure of the deterministic-equivalent CCO formulation. Although these approximations improve computational tractability, we must quantify their effect on the resulting network operating points. To this end, we perform ex-post validation using the solutions from the 25 simulations for each season. For every line and scheduling period, the squared-current term omitted from the optimization model is reconstructed as $\hat{\ell}{i,j,t}=(p_{i,j,t}^{2}+q_{i,j,t}^{2})/v_{i,t}$. This quantity is then used to estimate the neglected active-power losses, $R_{i,j}\hat{\ell}_{i,j,t}$, and the omitted nonlinear voltage-drop term, $(R_{i,j}^{2}+X_{i,j}^{2})\hat{\ell}_{i,j,t}$. In addition, we evaluate the optimized line flows against the exact circular thermal limit, quantifying both the maximum loading excess and the share of line--hour observations exceeding the exact limit. Table~\ref{tab:linearization_validation} reports the corresponding seasonal averages.

\begin{table}[h]
\centering
\caption{Ex-post validation of the lossless linearized network representation used in the CCO formulation.}
\label{tab:linearization_validation}
\renewcommand{\arraystretch}{1.15}
\setlength{\tabcolsep}{4pt}
\begin{tabular}{lccccc}
\toprule
\textbf{Season} & \textbf{Losses} & \textbf{Max. $\Delta v^{2}$} & \textbf{Mean loading} & \textbf{Max. loading excess} & \textbf{Violation share} \\
& [kWh] & [p.u.$^{2}$] & [\%] & [\%] & [\%] \\
\midrule
Autumn  & 0.0141 & $3.53\times10^{-8}$ & 38.3 & 0.0068 & 0.0244 \\
Spring  & 0.0139 & $3.53\times10^{-8}$ & 37.8 & 0.0076 & 0.0220 \\
Summer  & 0.0144 & $3.53\times10^{-8}$ & 40.0 & 0.0095 & 0.0415 \\
Winter  & 0.0143 & $3.53\times10^{-8}$ & 38.2 & 0.0051 & 0.0154 \\
\midrule
\textbf{Average} & \textbf{0.0142} & $\mathbf{3.53\times10^{-8}}$ & \textbf{38.6} & \textbf{0.0073} & \textbf{0.0258} \\
\bottomrule
\end{tabular}
\end{table}

The results show that the lossless network approximation has a negligible effect for the studied LV feeder. The estimated neglected losses range from 0.0139 to 0.0144~kWh across seasons, with an overall average of 0.0142~kWh, while the maximum omitted nonlinear voltage-drop term is only $3.53\times10^{-8}$~p.u.$^{2}$. Regarding the thermal-limit approximation, mean line loading ranges from 37.8\% in spring to 40.0\% in summer. The exact circular-limit check identifies only very small exceedances: the maximum loading excess ranges from 0.0051\% in winter to 0.0095\% in summer, while the proportion of line--hour observations exceeding the exact limit remains below 0.042\% in every season. Thus, for the specific characteristics and operating conditions of this case study, the neglected losses and voltage-drop terms are negligible, while the thermal-limit approximation remains within approximately 0.01\% of the exact circular boundary. These results indicate that the adopted linearized DN representation adequately approximates the conditions analyzed within the deterministic-equivalent CCO framework.

\subsection{Operational impact of uncertainty propagation}
As discussed in Section \ref{Propagation}, forecast uncertainty originates locally at individual agents but propagates through the radial DN as corresponding net-power deviations aggregate along the electrical paths. Consequently, evaluating uncertainty only at its point of origin does not fully capture its operational effect, since the resulting uncertainty margins may become more relevant in network elements that concentrate contributions from multiple downstream agents. The following analysis therefore examines how this propagated uncertainty is reflected in the CCO solution, first through its effect on line utilization and then through the community's coordinated net-power and storage response.

\begin{figure}[h]
    \centering
    \includegraphics[width=1.0\linewidth]{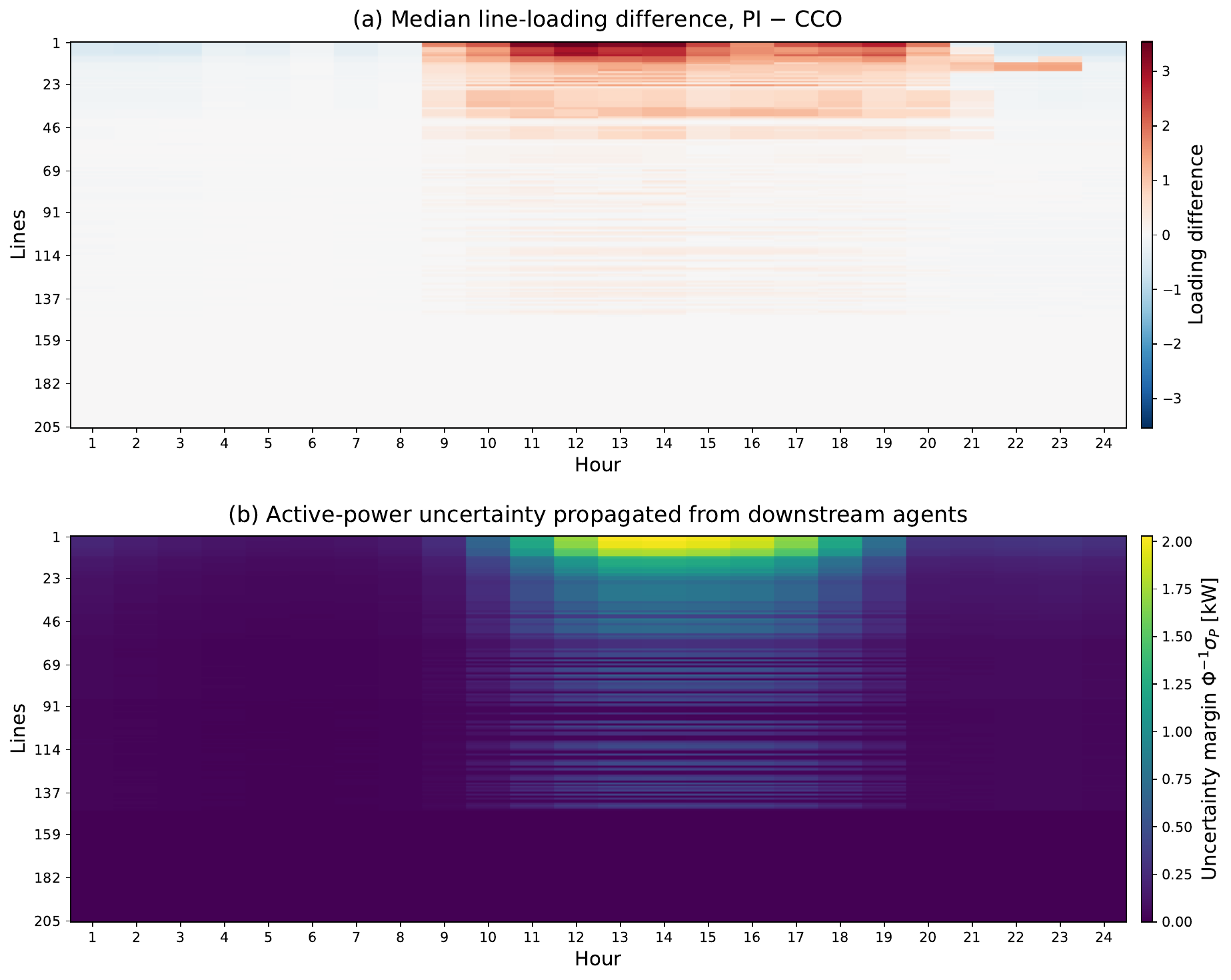}
    \caption{Line-loading differences between PI and CCO and propagated active power uncertainty across the DN.}
    \label{fig:heatmaps}
\end{figure}

Fig.~\ref{fig:heatmaps} compares, for each network line and hour, the difference between the PI and CCO loading levels with the corresponding propagated active-power uncertainty margin. In particular, panel (a) shows the median line-loading difference, defined as $PI_{ij,t}-CCO_{ij,t}$, where positive values (in red) indicate that the PI solution operates the corresponding line at a higher loading than the CCO solution, whereas negative values (in blue) indicate the opposite. Panel (b) shows the propagated active-power uncertainty margin, $\Phi^{-1}\sigma^{P}_{ij,t}$, resulting from aggregating the contribution of downstream forecast uncertainties along the radial network, so that upstream lines accumulate the uncertainty associated with a larger number of downstream agents.

The two heatmaps in Fig.~\ref{fig:heatmaps} show a clear temporal and spatial correspondence. The largest uncertainty margins occur mainly between hours 10 and 19, coinciding with the period of higher PV generation, and are concentrated in the upstream lines where downstream uncertainty is more strongly aggregated. During the same period, panel (a) shows predominantly positive loading differences in these lines, as indicated by the red regions, meaning that their loading is higher under PI than under CCO. This behavior results from the chance-constrained line-flow limits, which reserve part of the available line capacity to accommodate uncertainty. As $\sigma^{P}_{ij,t}$ increases, a larger share of the line capacity must remain available to accommodate potential forecast deviations, reducing the loading that the CCO solution can schedule. Although the propagated uncertainty originates from both demand and PV forecast errors, the temporal pattern observed is strongly concentrated during the solar period. To clarify the relative contribution of both uncertainty sources and their effect on community operation, the following analysis examines the evolution of the net-power and SOC uncertainty margins, $\sigma^{\Delta P}$ and $\sigma^{SOC}$.

\begin{figure}[h]
    \centering
    \includegraphics[width=1.0\linewidth]{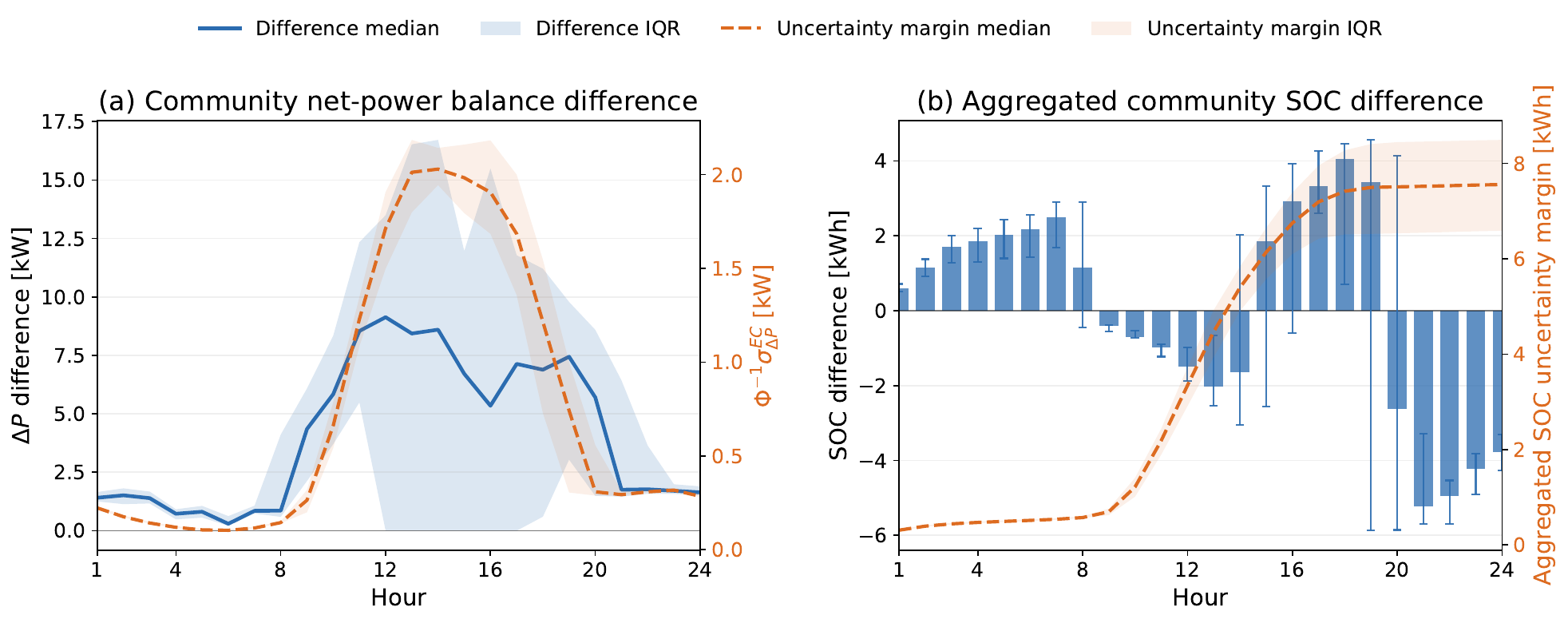}
    \caption{Community net-power and aggregated SOC differences between PI and CCO together with the corresponding uncertainty margins.}
    \label{fig:deltap_soc}
\end{figure}

Thus, Fig.~\ref{fig:deltap_soc}(a) compares the difference between the aggregated community net-power balance obtained under PI and the CCO solution, defined as $\Delta p^{PI}-\Delta ^{CCO}$, with the corresponding CCO uncertainty margin. Positive values therefore indicate that the PI solution presents a higher community net injection, or equivalently a lower net demand, than the CCO solution. During the first hours of the day, when PV production is negligible, the difference remains close to zero. The separation starts to increase after the beginning of the PV generation period and reaches its largest values during the central hours of the day. This evolution follows the increase in the uncertainty margin $\Phi^{-1}\sigma{\Delta P}^{EC}$, whose hourly profile is mainly driven by the PV uncertainty component. Demand uncertainty contributes throughout the entire horizon, while the additional increase during the solar period is associated with the PV contribution to $\sigma{\Delta P}$. As a result, the CCO adopts a lower net-power injection than the PI solution during hours when uncertainty in the community balance is higher. This behavior is consistent with the line-loading results discussed previously, where the propagated uncertainty during the same period leads the CCO to operate the upstream lines with a lower loading. 

Fig.~\ref{fig:deltap_soc}(b) shows the corresponding difference in the aggregated community SoC. A positive bar means that $soc^{PI}>soc^{CCO}$, indicating that the CCO has accumulated less energy up to that hour, either because it has charged less or discharged more than the PI solution. Conversely, a negative bar indicates that the CCO retains more stored energy than PI, corresponding to a higher cumulative charge or a lower cumulative discharge. The orange curve represents the propagated SOC uncertainty margin. Unlike net-power uncertainty, $\sigma_{soc}$ is cumulative by construction, since the uncertainty associated with PV generation and demand propagates through the storage balance over all previous periods. Its progressive increase during the day is therefore expected and reflects a gradual contraction of the admissible SOC operating range. The SOC difference exhibits four distinct operating periods that reflect how the CCO redistributes the available storage flexibility under the imposed uncertainty margins.

During the first hours, the positive difference indicates that the CCO maintains a slightly lower SOC than PI, creating additional storage capacity before the PV generation period. Around hours 9--14, the sign becomes negative as the CCO stores more energy than PI while PV generation and the associated net-power uncertainty increase. This behavior is consistent with the lower net injection observed in Fig.~\ref{fig:deltap_soc}(a), as part of the available PV energy is retained in the batteries rather than used or exported immediately. Between approximately hours 15 and 19, the difference becomes positive again, indicating that the CCO starts using the energy previously stored during the solar period as PV generation decreases. After hour 20, the sign reverses once more, and the CCO retains a higher SOC than PI. At this stage, the accumulated SOC uncertainty margin limits how much of the remaining stored energy can be used, whereas the PI solution, with the uncertainty already revealed, can exploit the available storage capacity more freely. The resulting sequence therefore reflects the more restrictive operation imposed by the CCO: storage absorbs part of the uncertainty during the PV period and redistributes energy over time, but SOC margins reduce its effective flexibility. Since the BESS provides the main source of temporal flexibility in the community, this restricted use of storage also helps explain why the CCO reaches a higher objective value than the PI benchmark.

\subsection{Sensitivity of CCO Performance to Uncertainty Parameterization}

The previous results evaluate the proposed CCO under a fixed uncertainty parameterization and show how the resulting uncertainty margins affect community and network operation. However, the behavior of a chance-constrained formulation depends on two key parameters: the forecast-error level $\epsilon$, which determines the magnitude of the standard deviations used to characterize PV-generation and demand uncertainty, and the confidence level $(1-\eta)$, which defines the probability with which the corresponding operational constraints must be satisfied. Assessing the sensitivity of the proposed formulation to these two parameters is therefore important to characterize the trade-off between economic performance and protection against forecast deviations, and to determine how progressively larger forecast errors and higher confidence requirements modify the resulting operating schedule.

To investigate this effect, we ran 15 simulations, each comprising one operating day for each of the four seasonal conditions. Within every simulation-season case, the same demand, PV, and electricity-price profiles were retained while the CCO was evaluated under 15 parameterizations obtained by combining three forecast-error levels, $\epsilon={5,10,15}\%$, with five confidence levels, $(1-\eta)={50,75,90,95,99}\%$, corresponding to $\Phi^{-1}(1-\eta)={0,0.674,1.282,1.645,2.326}$. This results in 900 CCO optimization runs in total. For each forecast-error level, we used the corresponding calibrated residual standard deviations described in Section~\ref{subsec:cco}, and scaled the resulting uncertainty margins by $\Phi^{-1}(1-\eta)$. The TSSP benchmark, which is independent of both $\epsilon$ and $(1-\eta)$, was solved once for each simulation-season case, resulting in 60 additional TSSP runs.

Table~\ref{tab:sensitivity} summarizes the average results obtained across the complete set of simulations and seasonal conditions. The table is organized into four main groups of indicators. The first group reports the design parameters varied in the sensitivity analysis, namely the forecast-error level $\epsilon$, the confidence threshold $(1-\eta)$, and the corresponding Gaussian quantile $\Phi^{-1}(1-\eta)$. The second group presents the objective function value, decomposed into the cost of energy purchased from the main grid and the revenue from energy sold to the grid. The third group reports the corresponding energy quantities, distinguishing between energy purchased from and sold to the external grid and the aggregate energy shared internally among community members. Finally, the last column reports the computational time required to solve each parameterization.

\begin{table*}[h]
\centering
\caption{Sensitivity analysis of the CCO framework: objective function decomposition, 
key operational indicators, and computational performance for parameterizations indexed by 
forecast-error level $\epsilon$ and confidence threshold $(1-\eta)$.}
\label{tab:sensitivity}
\renewcommand{\arraystretch}{1.18}

\begin{tabular}{
  c @{\hspace{7pt}} c @{\hspace{7pt}} r
  @{\hspace{14pt}}
  S[table-format=-2.2] @{\hspace{7pt}}
  S[table-format=2.2]  @{\hspace{7pt}}
  S[table-format=2.2]
  @{\hspace{14pt}}
  S[table-format=1.3] @{\hspace{7pt}}
  S[table-format=1.3] @{\hspace{7pt}}
  S[table-format=1.3]
  @{\hspace{14pt}}
  S[table-format=2.2]
}
\toprule
\multicolumn{3}{c}{\textbf{Design parameters}}
  & \multicolumn{3}{c}{\textbf{Objective decomposition [EUR]}}
  & \multicolumn{3}{c}{\textbf{Energy [MWh]}}
  & {\textbf{Time [s]}} \\
\cmidrule(r){1-3}
\cmidrule(lr){4-6}
\cmidrule(lr){7-9}
\cmidrule(l){10-10}

{$\epsilon$ [\%]}
  & {$(1{-}\eta)$ [\%]}
  & {$\Phi^{-1}(1{-}\eta)$}
  & {$F^{\mathrm{obj}}$}
  & {$C^{bg}$}
  & {$C^{sg}$}
  & { $E^{bg}$}
  & {$E^{sg}$}
  & {$E^{sh}$}
  & {Solve} \\
\midrule

\multirow{5}{*}{5}
  & 50 & 0.000 & -2.95 & 16.09 & 21.54 & 0.109 & 0.703 & 0.105 & 19.27 \\
  & 75 & 0.674 & -0.31 & 16.79 & 19.03 & 0.114 & 0.660 & 0.111 & 18.98 \\
  & 90 & 1.282 &  1.81 & 17.45 & 17.42 & 0.120 & 0.612 & 0.115 & 19.15 \\
  & 95 & 1.645 &  3.82 & 17.95 & 16.51 & 0.124 & 0.586 & 0.118 & 19.46 \\
  & 99 & 2.326 &  5.39 & 18.89 & 15.37 & 0.132 & 0.551 & 0.123 & 21.24 \\

\midrule

\multirow{5}{*}{10}
  & 50 & 0.000 & -2.95 & 16.09 & 21.54 & 0.109 & 0.703 & 0.105 & 19.09 \\
  & 75 & 0.674 &  2.24 & 17.56 & 17.31 & 0.121 & 0.609 & 0.116 & 19.31 \\
  & 90 & 1.282 &  7.03 & 19.25 & 14.33 & 0.135 & 0.526 & 0.124 & 19.53 \\
  & 95 & 1.645 &  8.80 & 20.45 & 13.07 & 0.142 & 0.494 & 0.129 & 21.48 \\
  & 99 & 2.326 & 13.18 & 24.01 &  9.45 & 0.160 & 0.450 & 0.138 & 19.17 \\

\midrule

\multirow{5}{*}{15}
  & 50 & 0.000 & -2.95 & 16.09 & 21.54 & 0.109 & 0.703 & 0.105 & 21.56 \\
  & 75 & 0.674 &  5.37 & 18.43 & 15.72 & 0.128 & 0.565 & 0.121 & 19.61 \\
  & 90 & 1.282 & 10.71 & 21.82 & 12.04 & 0.147 & 0.464 & 0.134 & 21.84 \\
  & 95 & 1.645 & 13.61 & 24.84 &  8.85 & 0.164 & 0.439 & 0.143 & 19.18 \\
  & 99 & 2.326 & 18.70 & 31.64 &  5.98 & 0.203 & 0.365 & 0.152 & 21.69 \\

\midrule

\multicolumn{3}{l}{\textit{TSSP}}
  & 3.05 & 19.38 & 15.37
  & 0.127 & 0.494 & 0.104 & 50.82 \\

\bottomrule
\end{tabular}

\medskip
\raggedright
\footnotesize
\textbf{Column definitions.}
$F^{\mathrm{obj}}$ = total objective function value;
$C^{bg}$ = grid procurement cost;
$C^{sg}$ = grid export revenue;
$E^{bg}$ = daily energy purchased from the main grid;
$E^{sg}$ = daily energy sold to the main grid;
$E^{sh}$ = aggregate intra-community shared energy.
\end{table*}

The results reveal a clear interaction between the assumed forecast-error level and the imposed confidence requirement. For a fixed confidence level, increasing $\epsilon$ progressively increases the objective function value, as larger uncertainty margins require a more conservative operating schedule. Similarly, for a fixed forecast-error level, increasing $(1-\eta)$ raises $\Phi^{-1}(1-\eta)$ and therefore enlarges the safety margins embedded in the chance constraints. This effect is reflected in higher grid-purchase costs, lower revenues from grid exports, and a progressive reduction in the CCO formulation's economic advantage relative to the TSSP benchmark. In particular, as the forecast-error level increases, the CCO solution departs progressively from the TSSP reference, while high confidence levels further reinforce this conservative behavior. Moreover, for forecast-error levels above $10\%$ combined with confidence levels above $90\%$, some individual simulations became infeasible. This occurs because the safety margins associated with the BESS state-of-charge constraints increase cumulatively over the scheduling horizon through $\sigma^{\mathrm{soc}}$, eventually leaving insufficient feasible operating range to satisfy the corresponding probabilistic SOC limits in some realizations.

\begin{figure}[h]
    \centering
    \includegraphics[width=1.0\linewidth]{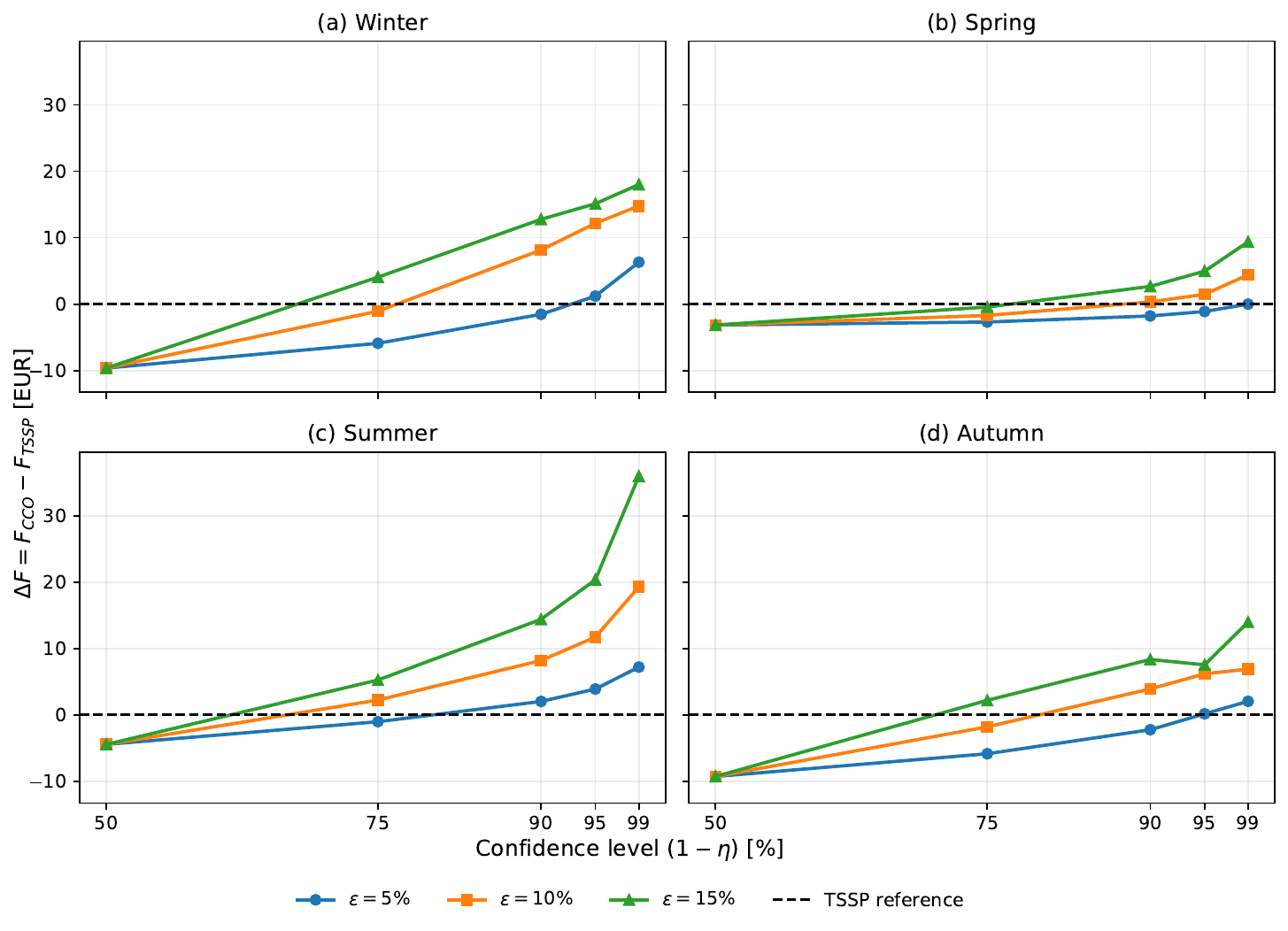}
    \caption{Seasonal comparison between CCO and TSSP in terms of the objective-function difference under different forecast-error and confidence levels.}
    \label{fig:CCO_Sensitivity}
\end{figure}

Additionally, when the results are disaggregated by season, the relative economic performance of the CCO and TSSP formulations is not uniform across operating conditions. Figure~\ref{fig:CCO_Sensitivity} illustrates this behavior by reporting, for each season, the difference between the objective function values obtained with both formulations, defined as $\Delta F = F_{\mathrm{CCO}}-F_{\mathrm{TSSP}}$. The horizontal axis represents the confidence threshold $(1-\eta)$, while the three curves correspond to forecast-error levels of $\epsilon=5\%$, $10\%$, and $15\%$. The horizontal dashed line at $\Delta F=0$ represents the TSSP reference: negative values indicate that the CCO achieves a lower objective function value than the TSSP, whereas positive values indicate a higher objective value and, therefore, a more conservative economic schedule relative to the stochastic benchmark.

The seasonal patterns shown in Fig.~\ref{fig:CCO_Sensitivity} are consistent with the relative behavior of the CCO and TSSP formulations discussed in Section~\ref{Comp_performance}. In particular, during summer, when the TSSP performed better economically, the CCO curves cross the TSSP reference at comparatively low confidence levels, especially for $\epsilon=10\%$ and $15\%$. This indicates that, under summer operating conditions, the CCO must operate with relatively low confidence requirements to preserve an economic advantage over the stochastic benchmark. Conversely, in autumn, where the TSSP showed larger deviations from the PI solution, the CCO remains below or close to the TSSP reference over a wider range of confidence levels, particularly for $\epsilon=5\%$. Winter and spring show more intermediate behavior, with spring showing comparatively smaller deviations between the two formulations over most of the analyzed parameter range.

More generally, the results indicate that the economic competitiveness of the CCO depends on the joint selection of the forecast-error level and the confidence requirement. For $\epsilon=5\%$, the proposed formulation remains economically competitive with the TSSP for confidence levels up to about $95\%$, since the objective-function differences remain negative or relatively small across seasons. In contrast, for larger forecast-error levels, particularly $\epsilon=10\%$ and $15\%$, increasing $(1-\eta)$ above $90\%$ systematically increases $\Delta F$, indicating progressively more conservative CCO schedules relative to the TSSP. This behavior becomes pronounced as the confidence level approaches $99\%$, where the combined effect of larger forecast uncertainty and more stringent probabilistic requirements substantially increases the economic cost of the CCO solution.

\subsection{Energy sharing}

Table~\ref{tab:stage4_energy_sharing} summarizes the community-level performance of the ex-post energy-sharing allocation. The reported indicators are computed independently for each simulation and season and are presented as the mean $\pm$ standard deviation across the 25 realizations. The self-sufficiency ratio (SSR) is calculated as $1-E^{bg}/E^{req}$, where $E^{req}$ is the aggregate member energy requirement, including demand and BESS charging. It therefore measures the fraction of these requirements supplied without relying on grid imports. In contrast, the self-consumption ratio (SCR) is calculated as $1-E^{sg}/E^{pool}$, where $E^{pool}$ is the collective PV generation and BESS discharge available for allocation. Hence, SSR is demand-oriented, whereas SCR is generation-oriented: a high SSR indicates low dependence on grid purchases, while a high SCR indicates that most of the available internal energy remains within the community instead of being exported. The sharing contribution (SC) measures the fraction of $E^{req}$ supplied through intra-community exchanges after applying the participation coefficients. Finally, the gross grid-exchange reduction (GER) compares the sum of member-level surpluses and deficits immediately after the $\beta_i$ allocation with the remaining grid imports and exports after internal sharing. It therefore quantifies the fraction of the gross external settlement that is avoided by matching post-allocation surpluses and deficits within the community.

\begin{table}[htbp]
\centering
\caption{Performance of the ex-post energy-sharing allocation. Values are reported as mean $\pm$ standard deviation across the simulated cases.}
\label{tab:stage4_energy_sharing}
\renewcommand{\arraystretch}{1.15}

\begin{tabular}{
    l
    @{\hspace{10pt}}
    S[table-format=2.2] @{\,$\pm$\,} S[table-format=2.2]
    @{\hspace{12pt}}
    S[table-format=2.2] @{\,$\pm$\,} S[table-format=2.2]
    @{\hspace{12pt}}
    S[table-format=2.2] @{\,$\pm$\,} S[table-format=2.2]
    @{\hspace{12pt}}
    S[table-format=2.2] @{\,$\pm$\,} S[table-format=2.2]
}
\toprule

\textbf{Season}
& \multicolumn{2}{c}{\textbf{SSR [\%]}}
& \multicolumn{2}{c}{\textbf{SCR [\%]}}
& \multicolumn{2}{c}{\textbf{SC [\%]}}
& \multicolumn{2}{c}{\textbf{GER [\%]}} \\

\midrule

Winter
& 70.46 & 4.27
& 64.02 & 20.14
& 14.61 & 5.80
& 30.85 & 16.72 \\

Spring
& 77.13 & 5.88
& 58.58 & 25.16
& 17.22 & 9.31
& 33.30 & 23.03 \\

Summer
& 79.77 & 2.48
& 40.94 & 17.16
& 9.30 & 6.44
& 14.13 & 16.94 \\

Autumn
& 71.78 & 7.33
& 53.39 & 25.71
& 11.40 & 5.98
& 21.58 & 18.69 \\

\midrule

Average
& 74.79 & 6.48
& 54.23 & 23.60
& 13.13 & 7.56
& 24.96 & 20.24 \\

\bottomrule
\end{tabular}

\medskip
\raggedright
\footnotesize
\textbf{Column definitions.}
SSR = self-sufficiency ratio;
SCR = self-consumption ratio;
SC = sharing contribution;
GER = gross grid-exchange reduction.
\end{table}

The seasonal results are consistent with the operating patterns identified in the previous analyses. Summer has the highest SSR, $79.77\%$, and the lowest variability ($\pm2.48\%$), reflecting the greater and more stable PV generation observed during this season. At the same time, summer exhibits the lowest SCR, $40.94\%$, showing that the larger internal energy pool exceeds the community's ability to absorb it and therefore leads to a greater fraction of exports. This also explains why SC and GER fall to $9.30\%$ and $14.13\%$, respectively: the dominant condition is an aggregate energy surplus rather than the coexistence of complementary member-level deficits and surpluses. Spring shows the opposite sharing behavior. Its SC reaches $17.22\%$, the highest seasonal value, while GER reaches $33.30\%$. Under these conditions, available internal generation is sufficient to create surpluses for part of the community but not enough to eliminate deficits for other members, creating more opportunities for internal matching after the $\beta_i$ allocation. Winter follows a similar pattern but with lower autonomy, whereas autumn combines an SSR of $71.78\%$ with the largest SSR dispersion ($\pm7.33\%$). This greater variability is consistent with the broader seasonal uncertainty previously observed in the demand, PV, and storage trajectories, for which the CCO remained closer to the PI benchmark than the TSSP formulation.
 
\section{Conclusions}
\label{sec:conclusions}
This paper developed a forecast-residual-based chance-constrained optimization framework for day-ahead LEC operation under PV-generation and demand uncertainty. The proposed formulation characterized uncertainty through the statistical properties of forecast residuals and embedded their effect into the scheduling problem through standard-deviation terms associated with local operational variables and their propagation across the DN. Under the adopted Gaussian assumption, the resulting chance constraints admitted an exact deterministic-equivalent reformulation, allowing the prescribed confidence level to be enforced through analytical safety margins while preserving the mixed-integer linear structure of the underlying scheduling problem. In addition, we incorporated an ex-post energy-allocation stage to separate the community's physical operation from the subsequent accounting of shared energy. Once the uncertainty-aware physical schedule was fixed, we aggregated available PV generation and BESS discharge into a common energy pool, allocated it among community members according to predefined participation coefficients, and then used it to identify and internally match post-allocation surpluses and deficits. 

We evaluated the proposed optimization framework on a reduced 206-bus European low-voltage DN with up to 55 community users. The computational experiments were designed to reproduce heterogeneous operating conditions through randomly generated demand and PV profiles, considering seasonal variability and different solar-generation conditions, including clear, partly cloudy, and cloudy days. The resulting synthetic forecast and realization samples were also used to calibrate the forecast-residual standard deviations required by the CCO and to construct the reference formulations, allowing evaluation of different uncertainty representations under consistent operating conditions. The numerical analysis showed that the CCO's main advantage emerged when forecast uncertainty remained moderate. At forecast-error levels near 5\%, the deterministic-equivalent formulation produced operating outcomes comparable to those from the TSSP benchmark while retaining computational requirements close to the deterministic formulation. Conversely, increasing either the forecast-error magnitude or the required confidence level enlarged these margins, progressively restricted the feasible operating region, and increased the conservativeness of the resulting schedules. The seasonal analysis further showed that the performance depended on how accurately the corresponding uncertainty representation captured the actual variability of demand and PV generation. Thus, the results highlighted that the effectiveness of a forecast-residual-based CCO depends fundamentally on the quality of the statistical characterization embedded in the $\sigma$ terms. The ex-post allocation stage showed that the proposed participation-coefficient mechanism can be applied without altering the physical schedule the CCO previously determined. After allocating the community energy pool, the method identified post-allocation surpluses and deficits and matched them internally, reducing the gross energy exchanges requiring settlement with the upstream grid.

Overall, the proposed approach provided a mechanism to translate the statistical information in forecast residuals directly into operational security margins, avoiding repeated construction and solution of scenario-expanded models while explicitly accounting for uncertainty propagation through DER and DN constraints. Its main limitation was therefore also linked to this statistical representation: the quality of the CCO solution depended on how accurately the adopted residual distributions and their standard deviations represented the underlying forecast uncertainty, together with the approximations introduced for uncertainty propagation through BESS and network variables. Future research should consequently focus on improving the estimation of these uncertainty parameters from increasing volumes of operational data. In particular, forecasting and machine-learning techniques could estimate time-, user-, and operating-condition-dependent forecast-error distributions and conditional standard deviations, allowing analytical safety margins to adapt more accurately to prevailing operating conditions. Alternative participation-coefficient rules could also be incorporated into the ex-post allocation stage without modifying the uncertainty-aware physical scheduling formulation.

\section*{CRediT Authorship Contribution Statement}

\noindent\textbf{Franco Cárdenas González:} Conceptualization, Methodology, Software, Validation, Formal analysis, Data curation, Writing---original draft, Visualization.

\noindent\textbf{Fernando García Muñoz:} Conceptualization, Methodology, Formal analysis, Writing---review and editing, Supervision, Project administration, Funding acquisition.

\section*{Funding}

This work has been supported by ANID FONDECYT Iniciación 11240745. The funders had no role in the study design, data collection and analysis, decision to publish, or preparation of the manuscript.

\bibliographystyle{elsarticle-num-names}
\bibliography{references}
\end{document}